\documentclass[11pt]{amsart}

\usepackage{pstricks,graphicx,pst-node}
\newtheorem{theorem}{Theorem}[section]
\newtheorem{lemma}[theorem]{Lemma}
\newtheorem{proposition}[theorem]{Proposition}
\newtheorem{corollary}[theorem]{Corollary}

\theoremstyle{definition}
\newtheorem{definition}[theorem]{Definition}

\theoremstyle{remark}
\newtheorem{remark}[theorem]{Remark}

\newenvironment{rmkT}[1]{\noindent\textit{Remarks on \reftab{#1}.}}{\vspace{\baselineskip}}

\newcommand{\reftab}[1]{Table \ref{tab:#1}}
\newcommand{\refeq}[1]{(\ref{eqn:#1})}

\newcommand{\refthm}[1]{Theorem \ref{thm:#1}}
\newcommand{\reflem}[1]{Lemma \ref{lem:#1}}
\newcommand{\refprop}[1]{Proposition \ref{prop:#1}}
\newcommand{\refcor}[1]{Corollary \ref{cor:#1}}

\newcommand{\refsec}[1]{Section \ref{sec:#1}}

\newcommand{\Sn}{\mathcal{S}_n}
\newcommand{\Ind}{\mathrm{Ind}}
\newcommand{\hG}{\widehat{G}}

\newcommand{\head}{\mathrm{h}}
\newcommand{\tail}{\mathrm{t}}

\newcommand{\SEP}{\mathrm{SEP}}
\newcommand{\linf}{l_{\infty}}

\newlength{\hsp}
\newlength{\vsp}
\newlength{\vspi}
\begin{document}


\title[A rule of thumb for riffle shuffling]{A rule of thumb for riffle shuffling}

\author[Assaf]{Sami Assaf}
\address{Department of Mathematics, Massachusetts Institute of Technology, 77 Massachusetts Avenue, Cambridge, MA 02139-4307}
\email{sassaf@math.mit.edu}

\author[Diaconis]{Persi Diaconis}
\address{Department of Statistics, Stanford University, 390 Serra Mall, Stanford, CA 94305-4065}

\author[Soundararajan]{K. Soundararajan}
\address{Department of Mathematics, Stanford University, 450 Serra Mall, Building 380, Stanford, CA 94305-2125}
\email{ksound@math.stanford.edu}

\subjclass[2000]{Primary 60B15; Secondary 60C05}

\date{\today}


\keywords{card shuffling, cutoff phenomenon, lumping of Markov chains,
  Poisson summation}

\begin{abstract}
  We study how many riffle shuffles are required to mix $n$ cards if
  only certain features of the deck are of interest, e.g. suits
  disregarded or only the colors of interest. For these features, the
  number of shuffles drops from $\frac{3}{2}\log_2 n$ to $\log_2
  n$. We derive closed formulae and an asymptotic `rule of thumb'
  formula which is remarkably accurate.
\end{abstract}

\maketitle

\section{Introduction}
\label{sec:intro}

In this paper we study the mixing properties of the
Gilbert-Shannon-Reeds model for riffle shuffling $n$
cards. Informally, the deck is cut into two piles by the binomial
distribution, and the cards are riffled together according to the
rule: if the left packet has $A$ cards and the right has $B$ cards,
drop the next card from the left packet with probability $A/(A+B)$
(and from the right packet with probability $B/(A+B)$). Continue until
all cards have been dropped. This defines a measure, denoted
$Q_2(\sigma)$, on the symmetric group $\Sn$. Repeated shuffles are
defined by {\em convolution powers}
\begin{equation}
  Q_2^{*k}(\sigma) = \sum_{\tau \in \Sn} Q_2(\tau)
  Q_2^{*(k-1)}(\sigma\tau^{-1}) .
\label{eqn:convolution}
\end{equation}
The {\em uniform distribution} is $U(\sigma) = 1/n!$. There are
several notions of the distance between $Q_2^{*k}$ and $U$: the {\em
  total variation distance}
\begin{equation}
  \| Q_2^{*k} - U \|_{TV} = \max_{A \subset \Sn} |Q_2^{*k}(A) - U(A)| =
  \frac{1}{2} \sum_{\sigma \in \Sn} |Q_2^{*k}(\sigma) - U(\sigma)| ,
\label{eqn:TVdist}
\end{equation}
and the {\em separation} and $\linf$ metrics
\begin{equation}
  \SEP(k) = \max_{\sigma} 1 - \frac{Q_2^{*k}(\sigma)}{U(\sigma)} \ ,
  \hspace{3em} \linf (k) = \max_{\sigma} \left| 1 -
    \frac{Q_2^{*k}(\sigma)}{U(\sigma)} \right| \ .
\label{eqn:SEPlinf}
\end{equation}

In widely cited works, Aldous \cite{AlDi1986} and Bayer and Diaconis
\cite{BaDi1992} show that $\frac{3}{2} \log_2(n) + c$ shuffles are
necessary and sufficient to make the total variation distance small,
while $2\log_2(n)+c$ shuffles are necessary and sufficient to make
separation and $\linf$ small.

The distances in \refeq{TVdist} and \refeq{SEPlinf} look at all
aspects of a permutation. In many card games, only some aspects of the
permutation matter. For example, in Black-Jack and Baccarat, suits are
irrelevant and all $10$'s and picture cards are equivalent; ESP card
guessing experiments use a Zener deck of $25$ cards with each of $5$
symbols repeated five times. It is natural to ask how many shuffles
are required in these situations. These questions are studied by
Conger and Viswanath \cite{CoVi2006,Conger2007,CoVi2007,CoVi}
who derive a framework, new formulae and remarkable numerical
procedures giving useful answers for cases of practical
interest. Their work is reviewed at the end of this introduction.

In this paper, we develop formulae and asymptotics for a deck of $n$
cards with $D_1$ cards labelled $1$, $D_2$ cards labelled $2$,
$\ldots$, $D_m$ cards labelled $m$. Most of the results are proved
from the deck starting `in order', i.e. with $1$'s on top through $m$'s
at the bottom. In \refsec{alternating}, we show that initial order can
change the conclusions.

In \refsec{Aspades}, we begin with $D_1=1$ and $D_2=n-1$. The
transition matrix for this case has interesting properties, rivaling
the `Amazing Matrix' in \cite{Holte1997}. Extending work of J.C. Reyes
\cite{Reyes2002}, we show that $\log_2 n
+c$ shuffles are necessary and sufficient for convergence in any of
our metrics.

\refsec{red-black} studies $D_1 = R$, $D_2 = B$, with, for example, $R
= B = 26$ modeling the red-black pattern for a standard $52$ card
deck. We derive a simple formula, first proved in \cite{Conger2007},
for $Q_2^{*k}(w)$ for any pattern $w$ and use this to again show that
$\log_2 n +c$ steps are necessary and sufficient for convergence to
uniformity. We find this surprising as following a single card
involves a state space of size $n$, reds and blacks involves a state
space of size $\binom{n}{n/2}$, and yet the same number of shuffles
are needed.

In \refsec{deck}, we treat the general case, deriving a formula which
can be used for some limited calculations. We also reprove a result of
Conger-Viswanath determining where the maximum for $\SEP$ and $\linf$
are achieved. A main result is a unified formula, our {\em rule of
  thumb}:

\begin{theorem}
  Consider a deck of $n$ cards with $D_i$ cards of type $i$, $1 \leq i
  \leq m$ with $D_i \geq d\ge 3$, $n = D_1+\cdots+D_m$. Then the separation
  distance after $k$ shuffles is
  \begin{displaymath}
    1 -(1+\eta) \frac{2^{k(m-1)}}{(n\!+\!1)\cdots(n\!+\!m\!-\!1)} \sum_{j=0}^{m-1}
    (-1)^{j} \binom{m-1}{j} \left( 1 - \frac{j}{2^k} \right)^{n+m-1},
  \end{displaymath}
  where $\eta$ is a real number satisfying 
  \[ 
  |\eta| \le \Big( 1+ \frac{n^2}{3(d-2)(2^k-m+1)^2}\Big)^{m-1}-1.
  \]
\label{thm:thumb}
\end{theorem}

This result does not depend on the individual details of the $D_i$ and
shows that the same number of shuffles are necessary and sufficient
for a variety of questions.  For numerical approximation, we set $\eta
=0$ and simply compute the single sum.  The bound on $\eta$ gives
explicit error estimates.  We demonstrate that the rule of thumb is
accurate for both single card and red-black problems studied in
earlier sections. Some numerical results are summarized below.

\newcommand{\rrc}{@{\extracolsep{10pt}}r}
 \begin{table}[ht]
  \begin{center}
    \caption{\label{tab:thumb} Rule of Thumb for the separation
      distance for $k$ shuffles of $52$ cards.}
    \begin{tabular}{@{\extracolsep{1pt}}c|
        @{\extracolsep{5pt}}r\rrc\rrc\rrc\rrc 
        \rrc\rrc\rrc\rrc\rrc 
        \rrc\rrc@{\extracolsep{0pt}}}
      \hline
      k & 1 & 2 & 3 & 4 & 5 & 6 & 7 & 8 & 9 & 10 & 11 & 12 \\ \hline
      \\ [-.5\vsp]
      BD-92 & 1.00 & 1.00 & 1.00 & 1.00 & 1.00 & 1.00 & 1.00 &
      .995 & .928 & .729 & .478 & .278 \\ [.5\vsp]
      blackjack & 1.00 & 1.00 & 1.00 & 1.00 & .999 & .970 & .834 &
      .596 & .366 & .204 & .108 & .056 \\ [.5\vsp]  
      $\clubsuit{\red\diamondsuit\heartsuit}\spadesuit$ & 1.00
      & 1.00 & .997 & .976 & .884 & .683 & .447 & .260 & .140 & .073
      & .037 & .019 \\ [.5\vsp]  
      {\red red}black & .962 & .925 & .849 & .708 & .508 & .317 & .179 &
      .095 & .049 & .025 & .013 & .006 \\
      \raisebox{-.5ex}{\includegraphics[height=3.3ex]{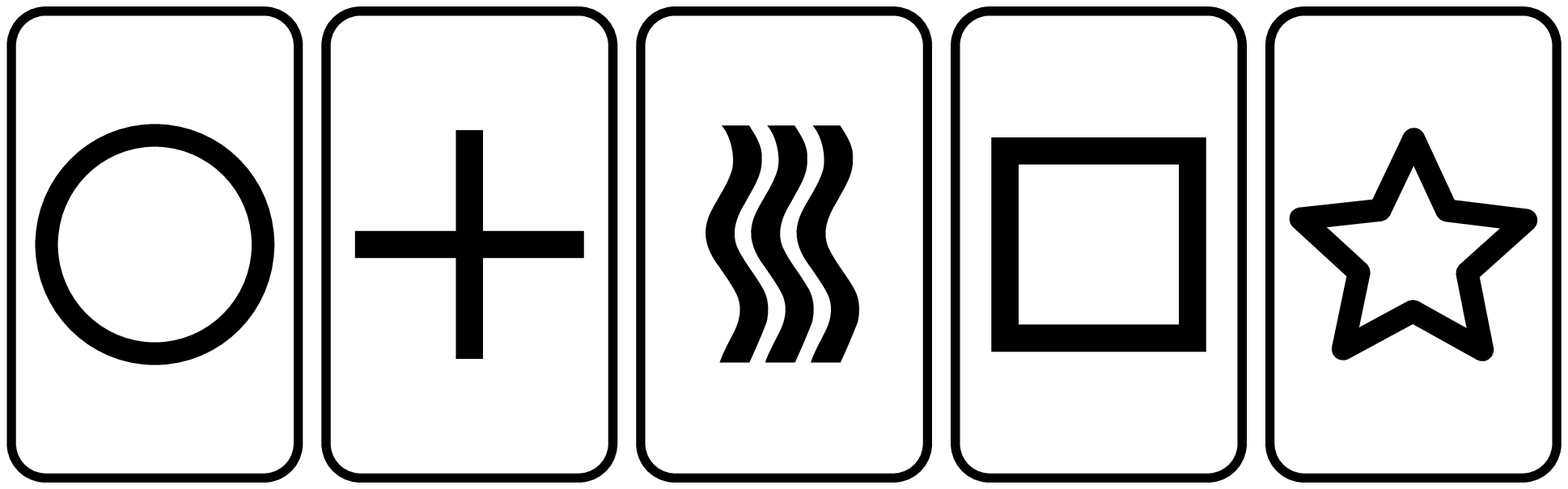}} & 
      1.00 & 1.00 & .993 & .943 & .778 & .536 & .321 & .177 & .093 &
      .048 & .024 & .012 \\ \hline 
    \end{tabular}
  \end{center}
\end{table}

\begin{rmkT}{thumb}
  The first row gives exact results from the Bayer-Diaconis formula
  for the full permutation group. The other numbers are from the rule
  of thumb. The single card or red-black numbers show that $6$
  shuffles achieve the the same separation as $12$ shuffles for the
  full deck. The Black-Jack (equivalently Baccarat) numbers suggest a
  savings of two or three shuffles, and the suit numbers lie in
  between. The final row is the rule of thumb for the Zener deck with
  25 cards, 5 cards for each of 5 suits.
\end{rmkT}

As explained at the end of Section~\ref{sec:deck}, the proof of
Theorem~\ref{thm:thumb} results from an approximation to an $m$-fold
iterated sum. Direct evaluation of this sum was sometimes possible
though it became intractable for many cases of interest. A referee
points out that we show that the sum is a coefficient in a product of
polynomials. Fast multiplication of polynomials results in a useful
polynomial time algorithm. One may ask, if the numbers can be computed
exactly, why bother with limit theorems and approximations? One answer
comes from understanding; The red/black and single card configurations
have very similar behaviors. This is surprising and calls out for
explanation. The rule of thumb formula explains why many different
configurations require the same number of shuffles.

In an appendix, we show that the processes studied below are quotient
walks with respect to Young subgroups of $\Sn$. We show how
representation theory can be used to derive results for features of
the random transposition random walk.

\subsection*{Literature review of riffle shuffles}

The basic shuffling model was introduced by Gilbert and Claude Shannon
in an unpublished report \cite{Gilbert1955}. The model was
independently introduced and studied by Jim Reeds in unpublished work
\cite{Reeds1976}. The first rigorous results are by Aldous
\cite{Aldous1983} who showed that asymptotically $\frac{3}{2}
\log_2(n)$ shuffles are correct for total variation. Separation
distance is introduced in connection with stopping time arguments in
Aldous and Diaconis \cite{AlDi1986}. They show that $2\log_2 n +c$
steps are necessary and sufficient for separation convergence. The
cutoff phenomena is first noticed in this paper as well. Recent work
on the cutoff phenomenon is in \cite{AlDi1987,GiSu2002,DiFu2008}. Our
work below adds several new examples to the list of problems where the
cutoff can be explicitly determined.

A generalization to $a$-shuffles is introduced by Bayer-Diaconis in
\cite{BaDi1992}. Here the deck is cut into $a$ packets by a
multinomial distribution, and then cards are dropped from packets with
probability proportional to packet size. Letting $Q_a(\sigma)$ denote
this measure, they show
\begin{equation}
  Q_a * Q_b = Q_{ab}.
\label{eqn:conv}
\end{equation}
Thus it is enough to study a single $a$-shuffle. The main result of
their paper is the simple formula
\begin{equation}
  Q_a (\sigma) = \frac{\binom{n+a-r}{n}}{a^n},
\end{equation}
where $r = r(\sigma)$ is the number of rising sequences in $\sigma$
($r(\sigma) = d(\sigma^{-1})+1$ with $d$ the number of descents in
$\sigma$). This allows simple closed form expressions for a variety of
distances. 

A number of extensions and variations have since developed. We will
not survey these here (see \cite{Diaconis2003} for a thorough
treatment) but mention that features of permutations are shown to
achieve the correct limiting distribution in fewer shuffles. For
example, $\frac{5}{6}\log_2 n + c$ suffice for the longest increasing
subsequences \cite{Fulman2002}, $\log_2(n)$ for the descent structure
\cite{DiFu2008}, $k_n \rightarrow \infty$ arbitrarily slowly for the
cycle structure \cite{DMP1995} and a single shuffle suffices for the
longest cycle \cite{DMP1995}. A recent addition is the work of Chen
and Saloff-Coste \cite{ChSa2008} studying random combinations of
$a$-shuffles for randomly varying $a$.

Mark Conger and D. Viswanath study the same type of problems as we do,
where cards are identified under the action of a subgroup. In
\cite{CoVi2006}, they lay out the basic problems, develop a formalism
for calculations involving descent polynomials (a generalization of
Eulerian polynomials), and use these to derive a closed formula for
the chance of a given arrangement after an $a$-shuffle for decks
labelled $\{1,2,\ldots,h,x^n\}$. This includes both our single card
case and the full deck case. They show that the probability of an
arrangement is
\begin{equation}
  \frac{1}{a^{n+h}} \sum_{m=r-1}^{a-1} \binom{m-r+h}{h-1} (a-m-1)^l
  (a-m)^{n-l}, 
\label{eqn:CV}
\end{equation}
with $r$ the number of cards labelled $c$, $1 \leq c \leq h$, that are
not preceded by a card labelled $c-1$ and $l$ the number of cards
labeled $x$ that precede the card labeled $h$. This elegant expression
can be analyzed asymptotically using the analytic techniques of
Sections 2-5 below. Their main results pertain to red-black decks
where they derive equivalence relations on configurations that have
the same probability. They point out that starting with the reds on
top or reds alternating with blacks can lead to different
conclusions. In a preliminary version of \cite{CoVi}, they give an
earlier proof of the exact formula for red-black decks found in
Theorem~\ref{thm:RB} below and good asymptotic approximations to the
total variation distance for following a single card. While this takes
$\log_2(n + c)$, they also prove the surprising result that the
seemingly similar problem of randomizing the current top card takes
half this many shuffles.

In \cite{CoVi2007}, the authors use their earlier work on descent
polynomials to develop a fascinating Monte Carlo procedure for
approximating the total variation distance. Our exact and asymptotic
calculations overlap theirs in many places, and in every case we find
their numbers spot on. This leads us to accept their estimates for
problems of deck hands at bridge where we have not found a way to do
exact calculations. The algorithms in \cite{CoVi2006,CoVi2007} can be
used to give polynomial time procedures for exact calculation of the
numbers in Table~\ref{tab:BD}. The authors have also proved some
complexity results showing that exact computation is intractable for
some of these problems. They have used their algorithms to calculate
an exact version of Table~\ref{tab:BD} above. Our numbers are based on
our rule of thumb. Their exact numbers agree to the accuracy given
except that for $k=1,2$ in the red-black category they get $.8898...$
and $.8897...$ and our approximation gives $.962$ and $.925$,
respectively. We are impressed (and thankful) for both their accurate
algorithms and for the accuracy of the rule of thumb. It is also worth
reporting that their stochastic approximation works in reasonable time
to give useful approximations for the analog of Table~\ref{tab:BD} for
total variation distance.

The results derived here add to the result of Conger-Viswanath in the
following ways. First, we present some new formulae (e.g. the
transition matrix for single card mixing or the red-black formula)
which allow exact computations. Second, we derive asymptotic
approximations for a variety of cases. Third, we supplement these
formulae and approximations with our unifying `rule of thumb'.

We mention the broad extensions of riffle shuffling to random walks on
hyperplane arrangements due to Bidigare, Hanlon and Rockmore (see
\cite{Diaconis2003} for a survey). The process induced by observing
which chamber of a sub-arrangement contains the present state of the
original walk is still Markov.  Rates of convergence for these
sub-arrangement walks are in \cite{AtDi2009}.

\subsection*{A Note on Metrics}

This paper focuses on the separation and $\linf$ metrics with some
attention on the more widely used total variation metric. These
metrics are related and often give `about the same answer'
asymptotically. For example, \cite{AlDi1987}[Proposition 5.13], it is
shown that separation distance $s(l)$ and total variation distance
$d(l)$ to uniformity after $l$ steps for a transitive Markov chain
satisfy:
\[
d(2l) \leq s(2l) \leq \phi(d(l))
\]
for $l \geq 1$ provided $d(l) \leq 1/8$, with $\phi(x) =
(1-(1-2x^{1/2}) (1 - x^{1/2})^2 \rightarrow 4 x^{1/2}$ as $x
\rightarrow 0$.  This shows (roughly) that $d(l)$ is small if and only
if $s(2l)$ is small.

In discussing total variation, separation and $\linf$ distances we
have sometimes heard the comment that these are worst case measures
that take their maximizing values on at a specific set or point. If
this set or point is not a particularly relevant configuration, the
distance may lose its relevance as well. One \textit{can} construct
artificial examples where this argument has merit, but worst case
bounds are conservative and the underlying inequalities show that they
hold for \textit{all} configurations. We suspect that in the present
natural context, there are many close by configurations that are also
`off'. This suggests interesting research questions.

Of course, the non asymptotic results depend on the details of
the metric. We do not think that there is `one right metric'. Often
$\linf$ and separation are easier to bound. Sometimes, the tails of
the distribution and rare points neglected by total variation and weak
star metrics are what is of interest. Unbounded functions such as the
number of correct guesses when cards are turned up sequentially
require separate treatment. The closed formulae reported here can be
used for any of these tasks.

Several other metrics are in wide-spread use. These include the
Chi-square or $l(2)$ distance (useful for reversible chains and
eigenvalue arguments), Entropy distance and various weak-star metrics
such as the maximum distance between the probability of balls in some
metric or the Wasserstein distance. For discussion and comparison see
\cite{GiSu2002,DiFu2008}.

\section{Following a single card}
\label{sec:Aspades}

Suppose one notices that the ace of spades is on the bottom of a deck
of $n$ cards. How many shuffles does it take until this one card is
close to uniformly distributed on $\{1,2,\ldots,n\}$?  This problem
was studied by J.C. Reyes \cite{Reyes2002}[Chapters 3 and 5]. He
derived the eigenvalues given below and also gave a coupling argument
that shows that $\log_2(n) + c$ shuffles suffice for total variation
convergence. A different proof of this result is by Fulman
\cite{Fulman2004}[Cor. 3.9], who derives it as a consequence of his
work on combining shuffles and random cuts. An asymptotic expansion
for the total variation appears in Conger and Viswanath
\cite{CoVi2007}. This shows that the upper bound cannot be
improved. In this section we elaborate on these results by studying
the transition matrix, giving its eigen values and vectors and giving
matching upper and lower bounds for the $\linf$, separation and total
variation distance. As shown in an appendix, under repeated shuffles a
single card moves according to a Markov chain. We begin by writing
down the transition matrix.

\begin{proposition}
  Let $P_a(i,j)$ be the chance that the card at position $i$ moves
  to position $j$ after an $a$-shuffle. For $1\leq i,j\leq n$,
  $P_a(i,j)$ is given by
  \begin{displaymath}
    \frac{1}{a^n} \sum_{k=1}^{a} \sum_{r=l}^{u} \!
    \binom{j\!-\!1}{r} \! \binom{n\!-\!j}{i\!-\!r\!-\!1} k^r
    (a-k)^{j-1-r}(k-1)^{i-1-r}(a-k+1)^{(n-j)-(i-r-1)}
  \end{displaymath}
  where $r$ ranges from $l = \max(0,(i+j)-(n+1))$ to $u = \min(i-1,j-1)$.
\label{prop:paij}
\end{proposition}

\begin{proof}
  We calculate $Q_a(j,i)$, the chance that an inverse $a$-shuffle
  brings the card at position $j$ to position $i$. For this to occur,
  the card at position $j$ may be labelled by $k$, $1 \leq k \leq a$.
  Then $r$ cards above this card may be labelled from $1$ to $k$. All
  will appear before the card at position $j$ in $\binom{j-1}{r}$
  ways. The remaining cards above must labelled from $k+1$ to
  $a$. Here $0 \leq r \leq \min(j-1,i-1)$. Also if $m$ cards below
  position $j$ are labelled from $1$ to $k-1$, then $m+r = i-1, m<n-j$
  and so $r \geq (i+j) - (n+1)$. Finally, $i-1-r$ cards below position
  $j$ must be labelled from $1$ to $k-1$ in $\binom{n-j}{i-r-1}$ ways,
  and the remaining cards must be labelled from $k+1$ to $a$.
\end{proof}

For example, the $n\times n$ transition matrices for $n=2,3$ are given
below. 
\begin{displaymath}
  \begin{array}{c}
    \displaystyle{\frac{1}{2a} \left( \begin{array}{cc}
        a+1 & a-1 \\ a-1 & a+1
      \end{array} \right)} \\[\vsp]
    \displaystyle{\frac{1}{6a^2} \left( \begin{array}{ccc}
        (a+1)(2a+1) & 2(a^2-1) & (a-1)(2a-1) \\
        2(a^2-1)    & 2(a^2+2) & 2(a^2-1) \\
        (a-1)(2a-1) & 2(a^2-1) & (a+1)(2a+1)
      \end{array} \right)}
  \end{array}
\end{displaymath}
Two other special cases to note are the extreme cases when $i=1$ or
$i=n$, which are given by
\begin{displaymath}
  P_a(1,j) = \frac{1}{a^n} \sum_{k=1}^{a} (a-k)^{j-1}(a-k+1)^{n-j},
  \hspace{3ex} 
  P_a(n,j) = \frac{1}{a^n} \sum_{k=1}^{a} k^{j-1}(k-1)^{n-j}.
\end{displaymath}

These single card transition matrices are studied by Ciucu
\cite{Ciucu1998} who gives a closed form for all $n$ when $a=2$:
\begin{displaymath}
  P_2(i,j) \ = \ \left\{ \begin{array}{rl}
      \frac{1}{2^n} \left(2^{i-1} + 2^{n-i}\right) & \mbox{if $i=j$}, \\[.5\vsp]
      \frac{1}{2^{n-j+1}} \binom{n-j}{i-1} & \mbox{if $i>j$}, \\[.5\vsp]
      P_2(n\!-\!i\!+\!1,n\!-\!j\!+\!1) & \mbox{if $i<j$}.
    \end{array} \right.
\end{displaymath}

These matrices share many properties of the `amazing matrix' developed
by Holte \cite{Holte1997}. See also \cite{DiFu2008} for connections
between Holte's amazing matrices and card shuffling. The following
Proposition is essentially due to Ciucu \cite{Ciucu1998}.

\begin{proposition} The transition matrices following a single card
  have the following properties:
  \begin{enumerate}
  \item they are {\em cross-symmetric}, i.e. $P_a(i,j) =
    P_a(n-i+1,n-j+1)$;
  \item $P_a \cdot P_b = P_{ab}$;
  \item the eigenvalues are $1, 1/a, 1/a^2, \ldots, 1/a^{n-1}$;
  \item the right eigen vectors are independent of $a$ and have the
    simple form: \\
    $V_m(i) = (i-1)^{i-1} \binom{m-1}{i-1} + (-1)^{n-i+m}
    \binom{m-1}{n-i}$ for $1/a^m$, $m \geq 1$.
  \end{enumerate}
\label{prop:amazing}  
\end{proposition}

\begin{proof}
  The cross-symmetry (1) follows from \refprop{paij}, and the
  multiplicative property (2) follows from the shuffling
  interpretation and equation \refeq{conv}. Property (1) implies that
  the eigen structure is quite constrained; see
  \cite{Weaver1985}. Properties (3) and (4) follow from results of
  Cuicu \cite{Ciucu1998}.
\end{proof}

\begin{remark}
  We note that Holte's matrix arose from studying the `carries
  process' of ordinary addition. Diaconis and Fulman \cite{DiFu2008}
  show that it is also the transition matrix for the number of
  descents in repeated $a$-shuffles. We have not been able to find a
  closer connection between the two matrices.
\end{remark}

From \refprop{paij} we obtain the following Corollary, which also 
follows as a special case of Theorem 2.2 in \cite{CoVi2006}.

\begin{corollary}
  Consider a deck of $n$ cards with the ace of spades starting at the
  bottom. Then the chance that the ace of spades is at position $j$
  from the top after an $a$-shuffle is
  \begin{equation}
    Q_a(j) \ = P_a(n,j) \ = \ \frac{1}{a^n} \sum_{k=1}^{a} (k-1)^{n-j} k^{j-1}.
    \label{eqn:Qai}
  \end{equation}
  \label{cor:Qai}
\end{corollary}

From the explicit formula, we are able to give exact numerical
calculations and sharp asymptotics for any of the distances to
uniformity. The results below show that $\log_2 n +c$ shuffles are
necessary and sufficient for both separation and total variation (and
there is a cutoff for these). This is surprising since, on the full
permutation group, separation requires $2\log_2 n+c$ steps whereas
total variation requires $\frac{3}{2}\log_2 n +c$. Of course, for any
specific $n$, these asymptotic results are just indicative.

\begin{table}[ht]
  \begin{center}
    \caption{\label{tab:BD}Distance to uniformity for a deck of $52$
      distinct cards.}
    \begin{tabular}{@{\extracolsep{1pt}}c | 
        \rrc\rrc\rrc\rrc\rrc \rrc\rrc\rrc\rrc\rrc \rrc\rrc@{\extracolsep{0pt}}}
      \hline
      & 1 & 2 & 3 & 4 & 5 & 6 & 7 & 8 & 9 & 10 & 11 & 12 \\ \hline \\
      [-.5\vsp]  
      $TV$ & 1.00 & 1.00 & 1.00 & 1.00 & .924 & .614 & .334 &
      .167 & .085 & .043 & .021 & .010 \\ [.5\vsp]
      $\SEP$ & 1.00 & 1.00 & 1.00 & 1.00 & 1.00 & 1.00 & 1.00 & .996 &
      .931 & .732 & .479 & .278 \\ [.5\vsp]  
      $\linf$ & $10^{53}$ & $10^{41}$ & $10^{29}$ &
      $10^{19}$ & $10^{12}$ & $10^{7}$ & $10^{5}$
      & 128 & 11.3 & 2.57 & .900 & .380 \\ \hline
    \end{tabular}
  \end{center}
\end{table}

\begin{table}[ht]
  \begin{center}
    \caption{\label{tab:AS}Distance to uniformity for a single card
      starting at the bottom of a $52$ card deck.}
    \begin{tabular}{@{\extracolsep{1pt}}c | 
        \rrc\rrc\rrc\rrc\rrc \rrc\rrc\rrc\rrc\rrc \rrc\rrc@{\extracolsep{0pt}}}
      \hline
       & 1 & 2 & 3 & 4 & 5 & 6 & 7 & 8 & 9 & 10 & 11 & 12 \\ \hline \\
       [-.5\vsp]  
      $TV$ & .873 & .752 & .577 & .367 & .200 & .103 & .052 &
      .026 & .013 & .007 & .003 & .002 \\ [.5\vsp]
      $\SEP$ & 1.00 & 1.00 & .993 & .875 & .605 & .353 & .190 &
      .098 & .050 & .025 & .013 & .006 \\ [.5\vsp]  
      $\linf$ & 25.0 & 12.0 & 5.51 & 2.37 & 1.02 & .460 & .217
      & .105 & .052 & .026 & .013 & .006 \\ \hline
    \end{tabular}
  \end{center}
\end{table}

\begin{table}[ht]
  \begin{center}
    \caption{\label{tab:AD}Distance to uniformity for a single card
      starting at the middle of a $52$ card deck.}
    \begin{tabular}{c | rrrr}
      \hline
       & 1 & 2 & 3 & 4 \\ \hline \\
       [-.5\vsp]  
      $TV$ & .494 & .152 & .001 & .000 \\ [.5\vsp]
      $\SEP$ & 1.00 & .487 & .003 & .000 \\ [.5\vsp]  
      $\linf$ & 1.92 & .487 & .003 & .000 \\ \hline
    \end{tabular}
  \end{center}
\end{table}

\begin{rmkT}{AS}
  We use \refprop{paij} to give exact results when $n=52$. For
  comparison, \reftab{BD} gives exact results for the full deck using
  \cite{BaDi1992}. Tables \ref{tab:AS} and \ref{tab:AD} show that it
  takes about half as many or fewer shuffles to achieve a given degree of
  mixing for a card at the bottom of the deck. For example, the widely
  cited `$7$ shuffles' for total variation drops this distance to
  $.334$ for the full ordering, but this requires only $4$ shuffles to
  achieve a similar degree of randomness for a single card at the
  bottom, and only $2$ for a single card starting in the
  middle. Similar statements hold for the separation and $\linf$
  metrics.
\end{rmkT}

For asymptotic results, we first derive an approximation to
separation. Since separation is an upper bound for total variation,
this gives an upper bound for total variation. Finally, we derive a
matching lower bound for total variation.

\begin{proposition}
  After an $a$-shuffle, the probability that the bottom card 
  is at position $i$ satisfies 
  \begin{displaymath}
    \frac{1}{a} \frac{\alpha^{n-i+1}}{1-\alpha^n} \le   Q_a(i) \le 
    \frac{1}{a} \frac{\alpha^{n-i}}{1-\alpha^{n-1}}, 
  \end{displaymath}
  where for brevity we have set $\alpha = 1-1/a$.  In particular, the 
  separation distance satisfies 
  \begin{displaymath}
    1- \frac{n}{a} \frac{\alpha^{n}}{1-\alpha^n} \ge \SEP(a) \ge
    1-\frac{n}{a} \frac{\alpha^{n-1}}{1-\alpha^{n-1}}.
  \end{displaymath}
  \label{prop:sepi-bound}
\end{proposition}

\begin{proof}  
  Since $k/(k-1) \ge a/(a-1)$ for all $1< k \le a$ we find that
  \begin{equation} 
    \alpha^{n-i} Q_a(n) \ge Q_a(i) \ge  \alpha^{-(i-1)} Q_a(1).
    \label{eqn:ineq}
  \end{equation}
  Therefore 
  \begin{displaymath}
    1 = \sum_{i} Q_a(i) \ge Q_a(1) \sum_{i=1}^{n} \alpha^{-(i-1)} =
    Q_a(1) a \alpha^{1-n} (1-\alpha^n), 
  \end{displaymath}
  so that
  \begin{displaymath}
    Q_a(1) \le \frac{1}{a} \frac{\alpha^{n-1}}{1-\alpha^n} \le \frac
    1a \frac{\alpha^{n-1}}{1-\alpha^{n-1}}. 
  \end{displaymath}
  Since $Q_a(n) =Q_a(1) +1/a$ it follows that $Q_a(n)\le \frac 1a
  \frac{1}{1-\alpha^{n-1}}$.  Using (\ref{eqn:ineq}) the desired upper
  bound for $Q_a(i)$ follows.
 
  Similarly,
  \begin{displaymath}
    1 =\sum_i Q_a(i) \le Q_a(n) \sum_{i=1}^{n} \alpha^{n-i} = Q_a(n)
    \frac{1-\alpha^n}{1-\alpha}, 
  \end{displaymath}
  so that 
  \begin{displaymath}
    Q_a(n) \ge \frac{1}{a} \frac{1}{1-\alpha^n}.
  \end{displaymath}
  Since $Q_a(1) =Q_a(n)-1/a$ it follows that $Q_a(1) \ge \frac{1}{a}
  \frac{\alpha^n}{1-\alpha^n}$, and from (\ref{eqn:ineq}) the desired
  lower bound for $Q_a(i)$ follows.  From (\ref{eqn:sepi}) and the
  above estimates we obtain our bounds on $\SEP(a)$.
\end{proof}

If $a = 2^{\log_2(n)+c} = n2^c$, then our result shows that the
$\SEP(a)$ is approximately
\begin{displaymath}
  1 -\frac{1}{2^c} \frac{e^{-2^{-c}}}{1-e^{-2^{-c}}}, 
\end{displaymath}
and for large $c$ this is $\approx 2^{-c-1}$.  The fit to the data in
\reftab{sep} is excellent: for example after ten shuffles of a
fifty-two card deck we have $2^{-c-1} = \frac{26}{1024}$ which is very
nearly the observed separation distance of $0.025$.

\begin{remark}
  \refprop{sepi-bound} gives a local limit for the probability that
  the original bottom card is at position $j$ \textit{from the
    bottom}. When the number of shuffles is $\log_2 n +c$, the density
  of this (with respect to the uniform measure) is asymptotically
  $z(c) e^{-j/2^c}$, with $z$ a normalizing constant ($z(c) =
  1/2^c(e^{j/2^c}-1)$). The result is uniform in $j$ for $c$ fixed,
  $n$ large.
\end{remark}

\begin{proposition} Consider a deck of $n$ cards with the ace of
  spades at the bottom.  With $\alpha= 1-1/a$, the total variation distance for the mixing
  of the ace of spades after an $a$-shuffle is at most
  \begin{displaymath}
    \frac{\alpha^{n+1}}{1-\alpha^n}
    -\frac{a\alpha^2(1-\alpha^{n-1})}{n(1-\alpha^n)}  +\frac{1}{n\log
      (1/\alpha)} \log \left(\frac an
      \frac{1-\alpha^n}{\alpha^{n+1}}\right), 
  \end{displaymath}
  and at least 
  \[ 
  \frac{\alpha^n}{1-\alpha^{n-1}} - \frac{a(1-\alpha^n)}{n\alpha(1-\alpha^{n-1})} 
  + \frac{1}{n\log (1/\alpha)} \log \Big( \frac{a}{n} \frac{1-\alpha^{n-1}}{\alpha^{n-1}}\Big).
  \]
\label{prop:TV-bound}
\end{proposition}

\begin{proof}  
  Let $Q_a(i)$ denote the probability that the ace of spades is at
  position $i$ from the top after an $a$ shuffle.  Note that $Q_a(i)$
  is monotone increasing in $i$, and let $i^*$ be such that $Q_a(i^*)
  < 1/n \le Q_a(i^*+1)$.  From Proposition \ref{prop:sepi-bound} we
  find that $i^*$ satisfies
  \begin{equation}
    \label{i*bound1}
    \frac{\alpha^{n-i^*+1}}{a(1-\alpha^n)} < \frac{1}{n} \le
    \frac{\alpha^{n-i^*-1}}{a(1-\alpha^{n-1})}, 
  \end{equation}
  so that 
  \begin{equation} 
    \label{i*bound2}
    \log \left(\frac an \frac{1-\alpha^{n-1}}{\alpha^{n-1}} \right)  
    \le i^* (\log 1/\alpha) \le \log \left( \frac{a}{n}
    \frac{1-\alpha^n}{\alpha^{n+1}}\right) 
  \end{equation} 
  
  From \refprop{sepi-bound} we have that the desired total variation
  is
  \begin{displaymath}
    \sum_{i\le i^*} \Big(\frac{1}{n} - Q_a(i)\Big) 
    \le \frac{i^*}{n} - \sum_{i\le i^*} \frac{\alpha^{n-i+1}}{a(1-\alpha^n)} 
    = \frac{i^*}{n} -\frac{\alpha^{n-i^*+1}}{1-\alpha^n} (1-\alpha^{i^*}),
  \end{displaymath}
  and also
  \[
  \sum_{i\le i^*} \Big(\frac{1}{n} - Q_a(i)\Big) 
    \ge \frac{i^*}{n} -\frac{\alpha^{n-i^*}}{1-\alpha^{n-1}} (1-\alpha^{i^*}).
    \]
  Using (\ref{i*bound1}) and (\ref{i*bound2}) we obtain the Proposition.
\end{proof}

\begin{remark}
  After $\log_2n +c$ shuffles, that is when $a=2^c n$,  \refprop{TV-bound}
  shows that the total variation distance is approximately (with $C=2^c$)
  \[ 
  C\log \Big( C (e^{1/C}-1)\Big) + \frac{1-C \log (e^{1/C}-1)}{(e^{1/C}-1)}.
  \]
  Thus when $c$ is `large and negative,' the total variation is close to $1$, 
  and when $c$ is large and positive, the total variation is close to $0$.  
  Thus  total variation and separation converge at the same rate. This is an
  asymptotic result and, for example, \reftab{AS} supports this.
\end{remark}

\begin{remark}
  From \refprop{sepw}, the $\linf$ distance is achieved for
  configurations with the ace of spades back on the
  bottom. \refprop{sepi-bound} gives a formula for this and the
  arguments for Propositions \ref{prop:sepi-bound} and
  \ref{prop:TV-bound} show that $\log_2 n + c$ shuffles are necessary
  and sufficient for convergence in $\linf$.
\end{remark}

\begin{remark}
 Similar, but more demanding, calculations show that if the ace of spades starts at
  position $i$, and $\max(i/n, (n-i)/n) \geq A > 0$ for some fixed
  positive $A$, then $\frac{1}{2} \log_2 n$ shuffles suffice for
  convergence in any of the metrics. We omit further details.
\end{remark}

\section{A red-black deck}
\label{sec:red-black}

We focus now on riffle shuffles of a deck consisting of $R$ red cards
and $B$ black cards. The purpose of this section is to give an
explicit description of $a$-shuffles of the deck with initial
configuration of red atop blacks. In Bayer-Diaconis \cite{BaDi1992},
the formula describing when an $a$-shuffle of $n$ distinct cards
results in a particular permutation has the simple form
\begin{displaymath}
  \frac{1}{a^n} \binom{a+n-r}{n},
\end{displaymath}
where $r$ is the number of rising sequences in the permutation. The
analysis for the red-black deck is markedly different. One indication
of this comes by noticing how likely the reverse deck is to occur. In
the case of permutations, the reverse deck has $n$ rising sequences,
and so the Bayer-Diaconis formula dictates that this configuration
cannot occur unless $a \geq n$. However, in the red-black case, the
reverse deck (blacks atop reds) may occur after a single $2$-shuffle
no matter the deck size. 

We begin by determining a formula for the chance of any arrangement
following an $a$-shuffle. This formula was proved earlier in
unpublished work of Conger and Viswanath \cite{CoVi}. We use the
result to derive the numbers in Table~\ref{tab:BD}. It also serves as
a simple case of the more complex argument in Section~\ref{sec:deck}
which also gives useful asymptotics.

\begin{theorem}
  Consider a deck with $R$ red cards on top of $B$ black cards. The
  probability that an $a$-shuffle will result in the deck
  configuration $w$ is
  \begin{equation}
    Q_a(w) = \frac{1}{a^{R+B}} \sum_{k=1}^{a} \sum_{j=1}^{R} 
    (k-1)^{R-j} k^{j-1} (a-k)^{b(j)} (a-k+1)^{B-b(j)}
    \label{eqn:Qa}
  \end{equation}
  where $b(j) =b_w(j)$ is the number of black cards above the $j$th red card
  in the deck $w$.
\label{thm:RB}
\end{theorem}

\begin{proof}
  The general formula for the probability of $w$ resulting from an
  $a$-shuffle is given by
  \begin{equation}
    \sum_{A_1 + \cdots + A_a = R+B}
    \frac{1}{a^n} \binom{R+B}{A_1,\ldots,A_a} \mathrm{prob}(w|A),
  \end{equation}
  where the sum is over all non-negative compositions $A =
  (A_1,A_2,\ldots,A_a)$ of $R+B$, i.e $A_i \geq 0$ and $A_1 + A_2 +
  \cdots + A_a = R+B$, and $\mathrm{prob}(w|A)$ denotes the
  probability that $w$ results from successively dropping cards from
  the piles $A_i$.  We break the sum into the following two cases:
  either there exists an integer $k$ such that $A_1 + A_2 + \cdots +
  A_k = R$ or not.

  Consider the case when the sum of the first $k$ piles is exactly
  $R$.  Then, the result of the subsequent riffle shuffle is equally
  likely to be any of the $\binom{R+B}{R}$ possible deck
  configurations. That is to say, given such a cut $A$,
  $\mathrm{prob}(w|A) = 1/\binom{R+B}{R}$ for every $w$. Therefore the
  contribution to $Q_a(w)$ from all such cuts is given by
  \begin{displaymath}
    \begin{array}{l}
      \displaystyle{
        \sum_{\substack{A_1+\cdots+A_a = R+B \\ \exists \ k \
            \mathrm{s.t.} \ A_1+\cdots+A_k = R}}
        \frac{1}{a^{R+B}} \binom{R+B}{A_1,\ldots,A_a} \frac{1}{\binom{R+B}{R}}}
      \\ [2\vsp]
      \hspace{3ex} = \displaystyle{\frac{1}{a^{R+B}} \sum_{k=1}^{a-1} 
        \sum_{A_k = 1}^{R} \! \sum_{\substack{A_{k+1}+\cdots+A_{a} = B \\
            A_1+\cdots+A_{k-1} = R - A_k}} \!\!\! \binom{R}{A_k} \!
        \binom{R-A_k}{A_1,\ldots,A_{k-1}} \! \binom{B}{A_{k+1},\ldots,A_{a}}}
      \\ [2\vsp]
      \hspace{3ex} = \displaystyle{\frac{1}{a^{R+B}}
        \sum_{k=1}^{a-1} (a-k)^B \sum_{A_k=1}^{R} \binom{R}{A_k} (k-1)^{R-A_k}}
      \\ [1.5\vsp]
      \hspace{3ex} = \displaystyle{\frac{1}{a^{R+B}}
        \sum_{k=1}^{a-1} (a-k)^B \left( k^R - (k-1)^R \right)}.
    \end{array}
  \end{displaymath}
  The choice to let $k$ be the {\em first} index such that $A_1 +
  \cdots + A_k = R$ is necessary in order to avoid over counting
  compositions with many $0$'s. This choice seemingly breaks the
  symmetry between $R$ and $B$ in the final formulation. However, the
  symmetric version may be obtained by taking $k$ to be the {\em last}
  index such that $A_1 + \cdots + A_k = R$. Finally, note that since
  $B \neq 0$, we may in fact take the sum over $k$ to range from $1$
  to $a$.

  Now consider the alternative case when there exists a pile
  (necessarily unique) containing both red and black cards. The
  assumption on $A$ amounts to the existence of integers $k,x,y$, with
  $1 \leq k \leq a$, $1 \leq x \leq R$, $1 \leq y \leq B$, such that
  $A_1 + \cdots + A_{k-1} = R-x$, $A_k = x+y$, and $A_{k+1} + \cdots +
  A_{a} = B-y$. Given such a cut $A$, $\mathrm{prob}(w|A) =
  r_{x,y}(w)/\binom{R+B}{R-x,x+y,B-y}$, where $r_{x,y}(w)$ denotes the
  number of rising subsequences consisting of $x$ red cards followed
  by $y$ black cards. The resulting contribution to $Q_a(w)$ from all
  such cuts is given by
  \begin{displaymath}
    \begin{array}{l}
      \displaystyle{\sum_{\substack{A_1+\cdots+A_a = R+B \\ \exists \ k \
            \mathrm{s.t.} \ A_1+\cdots+A_{k-1} < R \\ \mathrm{and} \
            A_{k+1}+\cdots+A_{a}<B}} 
        \frac{1}{a^{R+B}} \binom{R+B}{A_1,\ldots,A_a} \mathrm{prob}(w|A)}
      \\[2\vsp]
      \hspace{1ex} = \displaystyle{\frac{1}{a^{R+B}} \sum_{k=1}^{a}
        \sum_{x=1}^{R} \sum_{y=1}^{B} r_{x,y}(w) \hspace{-3ex}
        \sum_{\substack{A_1+\cdots+A_{k-1} = R-x \\
            A_{k+1}+\cdots+A_{a}=B-y}}  \hspace{-1ex} \binom{R-x}{A_1 \ldots
          A_{k-1}} \! \binom{B-y}{A_{k+1}\ldots A_{a}}} \\[1.5\vsp]
      \hspace{1ex} = \displaystyle{\frac{1}{a^{R+B}} \sum_{k=1}^{a}
        \sum_{x=1}^{R} \sum_{y=1}^{B} r_{x,y}(w) (k-1)^{R-x}
        (a-k)^{B-y}}.
    \end{array}
  \end{displaymath}
  For the final equation to make sense, we adopt the convention that
  $0^0=1$. 

  Let $b(j)$ denote the number of black cards above the $j$th red card
  in $w$. We may count rising subsequences of $w$ by the last red card
  used in the subsequence, giving the equation
  \begin{equation}  \label{eqn:rsubword}
    r_{x,y} (w) = \sum_{j=1}^{R} \binom{j-1}{x-1} \binom{B - b(j)}{y} .
    \end{equation}
  To see this, note that the first binomial coefficient counts the
  number choices of $x$ red cards before the $j$th red card, and the
  second binomial coefficient counts the number of choices for $y$
  black cards after the $j$th red card. Inserting this into the $x$
  and $y$ summations above gives
  \begin{displaymath}
    \begin{array}{l}
      \displaystyle{\frac{1}{a^{R+B}} \sum_{k=1}^{a} \sum_{x=1}^{R}
        \sum_{y=1}^{B} r_{x,y}(w) (k-1)^{R-x} (a-k)^{B-y}} \\
      \hspace{1ex} = \displaystyle{\frac{1}{a^{R+B}} \sum_{k=1}^{a}
        \sum_{j=1}^{R} \! \left(\sum_{x=0}^{R-1} \binom{j\!-\!1}{x}
          (k\!-\!1)^{R-x-1}\right) \!\! \left(\sum_{y=1}^{B}
          \binom{B\!-\!b(j)}{y} (a\!-\!k)^{B-y} \right)} \\ 
      \hspace{1ex} = \displaystyle{\frac{1}{a^{R+B}} \sum_{k=1}^{a}
        \sum_{j=1}^{R} (k\!-\!1)^{R-j}k^{j-1} 
          (a\!-\!k)^{b(j)} \left( (a\!-\!k\!+\!1)^{B-b(j)} - (a\!-\!k)^{B-b(j)}
          \right)}.
    \end{array}  
  \end{displaymath}

  The probability $Q_a(w)$ is obtained by adding the expressions in
  these two cases. Since 
  \begin{eqnarray*} 
    && \sum_{k=1}^{a} \sum_{j=1}^{R} (k-1)^{R-j} k^{j-1} (a-k)^B 
    = \sum_{k=1}^{a} k^{R-1} (a-k)^B \sum_{j=1}^{R}
    \left(\frac{k-1}{k}\right)^{R-j} \\  
    && \hspace{5ex} = \sum_{k=1}^{a} k^{R-1}(a-k)^B \frac{1-(1-1/k)^R}{1-(1-1/k)} 
    = \sum_{k=1}^{a} (a-k)^B (k^R-(k-1)^R), \\
  \end{eqnarray*} 
  we obtain the desired expression.
 \end{proof}

 Given \refeq{rsubword}, $Q_a$ gives a completely explicit description
 of $a$-shuffles, though this is difficult to evaluate for an
 arbitrary $w$. However, there are two special deck configurations for
 which $Q_a$ simplifies nicely, namely reds atop blacks (where
 $r_{x,y}(w) = \binom{R}{x}\binom{B}{y}$) and blacks atop reds (where
 $r_{x,y}(w) = 0$). By \refprop{sepw}, the formulae below can be used
 to give exact calculations for separation and $\linf$.

\begin{corollary}
  The probability of an $a$-shuffle resulting in the original deck
  configuration of reds atop blacks is
  \begin{displaymath}
    \frac{1}{a^{R+B}} \left(\sum_{k=1}^{a} \left( k^R - (k-1)^R
    \right) (a-k+1)^B  \right).
  \end{displaymath}
  The probability an $a$-shuffle resulting in the reverse deck
  configuration of blacks atop reds is
  \begin{displaymath}
    \frac{1}{a^{R+B}} \sum_{k=1}^{a-1} (a-k)^B \left( k^R - (k-1)^R
    \right) .
  \end{displaymath}
\label{cor:RB-BR}
\end{corollary}

Another special case to consider is tracking the position of a single
card starting at the bottom of the deck. For this case, taking $B=1$
and $R=n-1$ in \refeq{Qa} we recover  \refcor{Qai}.

Note that if instead we consider a single red card, i.e. $R=1$ and
$B=n-1$, starting at the top, then the distribution is the same. More
precisely, let $\widetilde{Q}_a(i)$ denote the chance that, say, the
$2$ of hearts is at position $i$ from the top of the deck after an
$a$-shuffle. Then it is easy to verify that $Q_a(i) =
\widetilde{Q}_a(n-i+1)$, which is just a special case of the
cross-symmetry in \refprop{amazing}.


Finally, consider the case of a single $2$-shuffle for an arbitrary
red-black deck. In this case, the left hand summand of \refeq{Qa}
reduces to a single term evaluating to $1$. For the right hand
summand, note that $k=1$ forces $x=R$, and $k=a$ forces $y=B$.

\begin{corollary}
  The probability of a $2$-shuffle resulting in a deck configuration
  $w$ is 
  \begin{equation}
    Q_2(w) = \frac{1}{2^{R+B}} \left( 2^{\head(w)} + 2^{\tail(w)} - 1 \right),
    \label{eqn:Q2}
  \end{equation}
  where $\head(w)$ denotes the number of red cards preceding the
  first black card in $w$, and $\tail(w)$ denotes the number of black
  cards following the final red card of $w$.
\label{cor:2shuffle}
\end{corollary}

Equation \refeq{Q2} can be used to give a simple formula for the
total variation after a single $2$-shuffle of a deck with $n$ red
cards and $n$ black cards. Here note that any two configurations with
the same number of red cards on top and black cards on bottom has the
same likelihood of occurrence. Therefore the total variation distance
after a single $2$-shuffle is given by
\begin{equation}
  \frac{1}{2} \left( \!
    \left(\frac{2^{n+1}-1}{2^{2n}} - \frac{1}{\binom{2n}{n}}\right)
    \! + \! \sum_{i=0}^{n-1} \sum_{j=0}^{n-1} \left| \frac{2^i + 2^j -
        1}{2^{2n}} - \frac{1}{\binom{2n}{n}} \right| \! \binom{2n \! -
      \! (i\!+\!j\!+\!2)}{n\!-\!(i\!+\!1)} \right)
\label{eqn:2sh}
\end{equation}
Using this formula, the total variation after a single $2$-shuffle of
a deck with $26$ red and $26$ black cards is $0.579$, which agrees
with the numerical approximations of Conger and Viswanath in
\cite{CoVi2006}. Conger and Viswanath have used their Monte Carlo
approximation to get useful total variation numbers. Their results
show that total variation convergence takes place much faster than
separation convergence in the red-black case. For 52 cards, after
$1,2,3,4$ shuffles it is $.579, .360, .208, .105$, respectively,
decreasing by a factor of two from then on.

Asymptotic results for the separation distance for red-black
configurations appear in the following section.

\section{Approach to uniformity in separation for general decks}
\label{sec:deck}

In this section we work with general decks containing $D_i$ cards
labelled $i$, $1 \leq i \leq m$. The following lemma shows that the
separation distance is always achieved by reversing the initial deck
configuration. Note this is equivalent to Theorem 2.1 from
\cite{CoVi2006}.

\begin{proposition}
  Let $D$ be a deck as above. After an $a$-shuffle of the deck with
  $1$'s on top down to $m$'s on bottom, the most likely deck
  configuration is this initial deck and the least likely
  configuration is the reverse deck $w^*$ with $m$'s on top down to
  $1$'s on the bottom. In particular, the separation distance is
  achieved for $w^*$.
\label{prop:sepw}
\end{proposition}

\begin{proof}
  Note first that the initial configuration can result from any
  possible cut of the deck into $a$ piles. Moreover, from any given
  cut of the deck, the identity is at least as likely to occur as any
  other configuration. The first assertion now follows. The only cuts
  of the initial deck which may result in $w^*$ are those containing
  no pile with distinct letters. However, for all such cuts, each
  rearrangement of the deck is equally likely to occur. Therefore
  $w^*$ minimizes $Q_a(w)$ and so maximizes $1 - Q_a(w)/U$.
\end{proof}

The explicit formula for $Q_a(w^{*})$ given in \refcor{RB-BR}
facilitates exact computations of $\SEP(a)$ for decks of practical
interest. Similarly, we can compute $Q_a(w^*)$ for an arbitrary deck
with $D_i$ $i$'s, $i=1,\ldots,m$.

\begin{theorem}
  Consider a deck with $n$ cards and $D_i$ cards labeled $i$,
  $i=1,\ldots,m$. Then the separation distance after an $a$-shuffle of
  the sorted deck ($1$'s followed by $2$'s, etc) is given by
  \begin{displaymath}
    1 \!-\! \frac{1}{a^n} \binom{n}{D_1\ldots D_m}
    \hspace{-2.5em} \sum_{\rule{0pt}{3ex} 0=k_0 < \cdots < k_{m-1} < a}
    \hspace{-2em} (a\!-\!k_{m\!-\!1})^{D_m} \! \prod_{j=1}^{m-1} 
    \left( (k_j \!-\! k_{j\!-\!1})^{D_j} \!-\! (k_j \!-\! k_{j\!-\!1}
      \!-\! 1)^{D_j} \right) .
  \end{displaymath}
\label{thm:sepG}
\end{theorem}

\begin{proof}
  From \refprop{sepw}, $w^*$ may only result from cuts with no pile
  containing distinct cards and any such cut is equally like to
  result in any deck. Therefore $Q_a(w^*)$ is given by
  \begin{displaymath}
    Q_a(w^*) = 
    \sum_{\substack{A_1 + \cdots + A_a = n \\ A \ \mathrm{refines} \ D}}
    \frac{1}{a^n} \binom{n}{A_1,\ldots,A_a}
    \frac{1}{\binom{n}{D_1,\ldots,D_m}},
  \end{displaymath}
  where `$A$ refines $D$' means there exist indices $k_1,\ldots,k_{m-1}$
  such that $A_1+\cdots+A_{k_1} = D_1$ and, for $i=2,\ldots,m-1$,
  $A_{k_{i-1} + 1} + \cdots + A_{k_i} = D_i$. Just as in the proof of
  \refthm{RB} we may take the $k_i$'s to be minimal so that the
  expression for $Q_a(w^*)$ simplifies to 
  \begin{equation}
    \frac{1}{a^n} \hspace{-1ex} \sum_{0=k_0 < \cdots <
      k_{m-1} < a} \hspace{-2em} (a\!-\!k_{m-1})^{D_m} \prod_{j=1}^{m-1}
    \left( (k_j \!-\! k_{j-1})^{D_j} - (k_j \!-\! k_{j-1} \!-\! 1)^{D_j} \right).
  \label{eqn:sepG}
  \end{equation}
  The result now follows from \refprop{sepw}.
\end{proof}

\begin{table}[ht]
  \begin{center}
    \caption{\label{tab:sep}Separation distance for $k$ shuffles of $52$ cards.}
    \begin{tabular}{@{\extracolsep{1pt}}c|
        @{\extracolsep{5pt}}r\rrc\rrc\rrc\rrc 
        \rrc\rrc\rrc\rrc\rrc 
        \rrc\rrc@{\extracolsep{0pt}}}
      \hline
      k & 1 & 2 & 3 & 4 & 5 & 6 & 7 & 8 & 9 & 10 & 11 & 12 \\ \hline
      \\ [-.5\vsp] 
      BD-92 & 1.00 & 1.00 & 1.00 & 1.00 & 1.00 & 1.00 & 1.00 & .995 &
      .928 & .729 & .478 & .278 \\ [.5\vsp] 
      blackjack & 1.00 & 1.00 & 1.00 & 1.00 & .999 & .970 &
      .834\makebox[0pt]{$\ ^*$} & .596\makebox[0pt]{$\ ^*$} &
      .366\makebox[0pt]{$\ ^*$} & .204\makebox[0pt]{$\ ^*$} &
      .108\makebox[0pt]{$\ ^*$} & .056\makebox[0pt]{$\ ^*$}\\ [.5\vsp]   
      $\clubsuit{\red\diamondsuit\heartsuit}\spadesuit$ & 1.00 & .997
      & .997 & .976 & .884 & .683 & .447 & .260 & .140 & .073 &
      .037\makebox[0pt]{$\ ^*$} & .019\makebox[0pt]{$\ ^*$} \\
      [.5\vsp]   
      A$\spadesuit$ & 1.00 & 1.00 & .993 & .875 & .605 & .353 & .190 &
      .098 & .050 & .025 & .013 & .006 \\ [.5\vsp]  
      {\red red}black & .890 & .890 & .849 & .708 & .508 & .317 & .179
      & .095 & .049 & .025 & .013 & .006 \\ [.3\vsp] 
      \raisebox{-.5ex}{\includegraphics[height=3.5ex]{zener.eps}} &
      1.00 & 1.00 & .993 & .943 & .778 & .536 & .321 & .177 &
      .093\makebox[0pt]{$\ ^*$} & .048\makebox[0pt]{$\ ^*$} &
      .024\makebox[0pt]{$\ ^*$} & .012\makebox[0pt]{$\ ^*$} \\
      \hline
    \end{tabular}
  \end{center}
\end{table}

\begin{rmkT}{sep}
  We calculate $\SEP$ after repeated $2$-shuffles for various decks
  using \refthm{sepG}: (blackjack) $9$ ranks, say A$23456789$, with
  $4$ cards each and another rank, say $10$, with $16$ cards;
  ($\clubsuit\diamondsuit\heartsuit\spadesuit$) $4$ distinct suits,
  say clubs, diamonds, hearts and spades, of $13$ cards each;
  (A$\spadesuit$)the ace of spades and $51$ other cards; (redblack) a
  two color deck with $26$ red and $26$ black cards; and
  (\raisebox{-.8ex}{\includegraphics[height=3ex]{zener.eps}}) a deck
  with 5 cards in each of 5 suits. The entries in \reftab{sep}
  indicated by $^*$ were provided by the referee using
  Remark~\ref{rmk:ROT} below.
\end{rmkT}

\refprop{sepw} may be used with the Conger-Viswanath formula in
\refeq{CV} to give a simple expression for separation after an
$a$-shuffle for a deck of size $h+n$ with cards labelled
$1,2,\ldots,h$ and $n$ cards labelled $x$:
\begin{displaymath}
  \SEP(a) = 1 - \frac{(n+h)\cdots(n+1)}{a^{n+h}} \sum_{k=h-1}^{a-1}
  \binom{k}{h-1} (a-1-k)^n.
\end{displaymath}

Now we derive a basic asymptotic tool, \refprop{Poissonsum}, which
allows asymptotic approximations for general decks. As motivation,
consider again the case of one card mixing, i.e. begin with $n$ cards
with the ace of spaces at the bottom of the initial deck. How many
shuffles are required to randomize the ace of spades?  Recall from
\refcor{Qai} that the chance that the ace of spades is at position $i$
from the top after an $a$-shuffle is given by
\begin{displaymath}
  Q_a(i) \ = \ \frac{1}{a^n} \sum_{k=1}^{a} (k-1)^{n-i} k ^{i-1},
\end{displaymath}
with the convention $0^0=1$. Therefore from \refprop{sepw}, we have
\begin{equation}
  \SEP(a) = 1 - nQ_a(1) = 1 - \frac{n}{a^n} \sum_{k=1}^{a}(k-1)^{n-1}.
  \label{eqn:sepi}
\end{equation}
Exact calculations when $n=52$ are given in \reftab{sep}. 

\begin{proposition}
  Let $a$ be a positive real number, and let $r$ and $s$ be natural
  numbers with $r$, $s \ge 2$.  Let $\xi$ be a real number in $[0,1]$.
  Then
  \begin{eqnarray*}
    S(a,\xi;r,s)&:=& \frac{1}{a^{r+s}} \sum_{0\le k \le a-\xi}
    (k+\xi)^{r} (a-k-\xi)^s  \\
    &=& a\frac{r! s!}{(r+s+1)!} + \frac{\theta}{6a} \frac{r!
      s!}{(r+s-1)!} \Big( \frac{1}{r-1}+\frac{1}{s-1}\Big),\\ 
  \end{eqnarray*}
  where $\theta$ is a real number in $[-1,1]$. 
\label{prop:Poissonsum}  
\end{proposition}

\begin{proof}  
  Put $f(x) = x^r (1-x)^s$ for $x\in [0,1]$ and $f(x)=0$ otherwise.
  The sum that we wish to evaluate is
  \begin{equation} 
    \label{eqn:PS}
    \sum_{k \in {\mathbb Z}} f((k+\xi)/a) = a \sum_{\ell \in {\mathbb
        Z} } {\hat f}(a\ell) e(\ell \xi), 
  \end{equation}
  by the Poisson summation formula.  Here, we write $e(x)=e^{2\pi i
    x}$ and ${\hat f}(y) = \int_{-\infty}^{\infty} f(x) e(-xy) dx$
  denotes the Fourier transform.

  Now note that 
  \begin{equation}
    \label{eqn:PS1}
    {\hat f}(0) = \int_0^1 x^r (1-x)^s dx = \frac{r! s!}{(r+s+1)!}.
  \end{equation}
  Further 
  \begin{eqnarray*}
    {\hat f}(y) & = & \int_0^1 x^r (1-x)^s e^{-2\pi i xy} dx =
    \frac{1}{2\pi i y} \int_0^1 f^{\prime}(x) e^{-2\pi i xy} dx \\
    & = & \frac{1}{(2\pi i y)^2} \int_0^1 f^{\prime \prime}(x) e^{-2\pi i xy}dx,
  \end{eqnarray*}
  upon integrating by parts twice, and since $r$, $s \ge 2$ we have
  $f(0)=f^{\prime}(0)=f(1) = f^{\prime}(1)=0$.  Therefore
  \[
  |{\hat f}(y)|\le \frac{1}{4\pi^2 y^2 }\int_0^1 |f^{\prime \prime}(x)|dx.
  \] 
  Now 
  \[
  f^{\prime\prime}(x) = \Big(\frac rx -\frac{s}{1-x}\Big)^2 x^r (1-x)^s - \Big(\frac r{x^2} 
  +\frac s{(1-x)^2}\Big) x^r (1-x)^s,
  \] 
  and so  
  \begin{displaymath}
    \begin{array}{l}
      \displaystyle{\int_0^1 |f^{\prime \prime}(x)| dx} \\
      \hspace{2ex} \displaystyle{\leq \int_0^1 \left( \frac rx
      -\frac {s}{1-x}\right)^2  x^r (1-x)^s dx + \int_0^1 \left( \frac
      r{x^2} +\frac{s}{(1-x)^2}\right) x^r (1-x)^s dx}  \\
      \hspace{2ex} \displaystyle{= \frac{r!s!}{(r+s-1)!} \left( \frac{2}{r-1}+
        \frac{2}{s-1}\right)}.
    \end{array}
  \end{displaymath}    
  Combining the above estimates with (\ref{eqn:PS}) and
  (\ref{eqn:PS1}) we conclude that our sum equals
  \begin{displaymath}
    a\frac{r!s!}{(r+s+1)!} + \frac{\theta}{2\pi^2 a}
    \frac{r!s!}{(r+s-1)!} \left(\frac{1}{r-1}+\frac{1}{s-1} \right) 
    \sum_{\substack{\ell \in {\mathbb Z} \\ \ell\neq 0}} \frac{1}{\ell^2} 
  \end{displaymath}
  for some $\theta \in [-1,1]$.  Since $\sum_{\ell=1}^{\infty}
  \ell^{-2}=\pi^2/6$ the Proposition follows.
\end{proof} 

Now suppose we have $n$ red cards and $n$ black cards, so $2n$ cards
altogether, with the red cards starting on top. In this case, the
uniform distribution $U(w) = U = 1/\binom{2n}{n}$. Again we use
\refprop{sepw} this time with \refcor{RB-BR} to give a formula for the
separation distance,
\begin{equation}
  \SEP(a) = 1 - \binom{2n}{n}Q_a(w^*) = 1 -
  \frac{\binom{2n}{n}}{a^{2n}} \sum_{k=1}^{a-1}(a-k)^{n} \left(k^n -
    (k-1)^n\right) 
\label{eqn:sepnn}
\end{equation}
For exact computations when $2n=52$, see \reftab{sep}.  We now use
\refprop{Poissonsum} to calculate asymptotic expressions for
this separation distance.

\begin{corollary}
  For $2n$ cards starting with $n$ red cards on top, we have, with
  $\alpha=1-1/a$
  \begin{displaymath}
    \SEP(a)=1- \frac{a}{2n+1} (1-\alpha^{2n+1}) + \frac{2\theta}{3a}
    \frac{n}{(n-2)} (1-\alpha^{2n-1}), 
  \end{displaymath}
  for some real number $\theta \in [-1,1]$.
  In particular, for $n$ large with $a = 2^{\log_2(2n)+c}$,
  \begin{displaymath}
    \SEP(a) = 1 - 2^{c} \left(1 - e^{-2^{-c}}\right) + O\Big( \frac{1}{a}\Big).
  \end{displaymath}
  \label{cor:sepnn-bound}
\end{corollary}

\begin{proof} 
  Note that 
  \begin{displaymath}
    \begin{array}{l}
      \displaystyle{\frac{1}{a^{2n}}\sum_{k=1}^{a} (a-k)^n  (k^n-(k-1)^n)} \\
      \hspace{2em} \displaystyle{= \frac{1}{a^{2n}} \sum_{k=1}^{a}
        (a-k)^n  \int_0^1 n(k-1+\xi)^{n-1} d\xi} \\ 
      \hspace{2em} \displaystyle{= \frac{n}{a^{2n}} \int_0^1
        \sum_{k=0}^{a-1} (a-1+\xi -(k-1+\xi))^n (k-1+\xi)^{n-1} d\xi} .
    \end{array}
  \end{displaymath}
  Using \refprop{Poissonsum} we see that the inner sum over $k$ above
  equals
  \[ 
  (a-1+\xi)^{2n} \frac{n!(n-1)!}{(2n)!} + (a-1+\xi)^{2n-2} \frac{\theta}{6} \frac{n!(n-1)!}{(2n-2)!} \Big(\frac{1}{n-1}+\frac{1}{n-2}\Big).
  \]
  Using these observations in (\ref{eqn:sepnn}) we obtain that
  $\SEP(a)$ is given by
  \[ 
  1- \int_0^1 \Big(\frac{a-1+\xi}{a}\Big)^{2n} d\xi + \frac{\theta}{6a^2}\frac{2n(2n-1)(2n-3)}{(n-1)(n-2)} 
  \int_0^1 \Big(\frac{a-1+\xi}{a}\Big)^{2n-2} d\xi. 
  \]
  With a little calculus the Corollary follows.
\end{proof}
 
The approximation 
\begin{equation}
  \binom{2n}{n} \sum_{k=1}^{a} (a-k)^n (k^n-(k-1)^n)  \approx
  \frac{a^{2n+1}-(a-1)^{2n+1}}{2n+1} 
  \label{approx}
\end{equation}
which is the basis of our Corollary above is more accurate than
suggested by the simple error bounds that we have given.  For example,
when $n=26$ and $a=16$, the actual separation distance (given in
\reftab{sep}) differs from the approximation of the Corollary by about
$7 \times 10^{-12}$.  Put differently, note that the LHS and the RHS
of (\ref{approx}) are both polynomials in $a$ of degree $2n$, and in
fact the coefficients of both polynomials match for all degrees
between $n$ and $2n$.
 
Before moving to general decks, we establish a generalization of 
\refprop{Poissonsum}. 
 
 \begin{proposition}  Let $m\ge 2$ and $a$ be natural numbers, let 
 $\xi_1$, $\ldots$, $\xi_m$ be real numbers in $[0,1]$.  
 Let $r_1$, $\ldots$, $r_m$ be natural numbers all 
 at least $r\ge 2$.  Let 
 \[
 S_m(a; \underline{\xi}, \underline{r}) =
 \sum_{\substack{a_1,\ldots,a_m \ge 0 \\ a_1+\ldots+a_m=a}}
 (a_1+\xi_1)^{r_1} \cdots (a_m+\xi_m)^{r_m}.
 \]
 Then 
 \begin{displaymath}
   \begin{array}{l}
     \displaystyle{\Big| S_{m}(a;\underline{\xi},\underline{r}) 
       -\frac{r_1! \cdots r_m!}{(r_1+\ldots+r_m+m-1)!}
       (a+\xi_1+\ldots+\xi_m)^{r_1+\ldots+r_m+m-1}\Big|} \\
     \displaystyle{\leq r_1! \cdots r_m! \sum_{j=1}^{m-1}
       \binom{m-1}{j} \Big(\frac{1}{3 (r-1)}\Big)^j  \frac{ (a
         +\xi_1+\ldots+\xi_m)^{r_1+\ldots+r_m
           +m-1-2j}}{(r_1+\ldots+r_m+m-1-2j)!}}.
   \end{array}
 \end{displaymath}
 \label{prop:GPS}
 \end{proposition} 
 \begin{proof}
   We establish this by induction on $m$.  The case $m=2$ follows from
   \refprop{Poissonsum}, taking there $a$ to be what we would now call
   $a+\xi_1+\xi_2$.  Let now $m\ge 3$ and suppose the result has been
   established for $m-1$ variables.  Now
 \begin{equation}
 S_m(a;\xi,r) = \sum_{a_1 =1}^{a+\xi_2+\ldots+\xi_m-1} a_1^{r_1} S_{m-1} (a-a_1;\underline{\tilde{\xi}},\underline{\tilde{r}}) 
 \label{eqn:ind}
 \end{equation}
 with $\underline{\tilde{\xi}}= (\xi_2, \ldots, \xi_m)$ and $\underline{\tilde{r}} = (r_2,\ldots,r_m)$, and 
 interpreting the terms with $a_1 \ge a$ as being $0$. 
 Using the induction hypothesis we have that 
 \begin{displaymath}
   \begin{array}{l}
   \displaystyle{\Big| S_{m-1}(a-a_1;\underline{\tilde\xi},
     \underline{\tilde r}) \ -} \\
   \hspace{2em}
   \displaystyle{%
     \frac{r_2! \cdots r_m!}{(r_2+\ldots+r_m+m-2)!}
     (a-a_1+\xi_2+\ldots+\xi_m)^{r_2+\ldots+r_m+m-2} \Big|} \\
   \displaystyle{
     \le r_2! \cdots r_m! \! \sum_{j=1}^{m-2} \!\! \binom{m\!-\!2}{j} \!
     \Big(\! \frac{1}{3(r\!-\!1)}\!\Big)^j
     \frac{(a\!-\!a_1\!+\!\xi_2\!+\! \ldots\! +\!\xi_m)^{r_2+\ldots+r_m+m-2-2j}}%
     {(r_2\!+\!\ldots\!+\!r_m\!+\!m-2-2j)!}}. 
 \end{array}
\end{displaymath}
Note that the above estimate is valid even if $a+\xi_2 +\ldots +\xi_m
-1 \ge a_1 \ge a$ since the RHS is larger than the main term that is
being subtracted in the LHS.  We use this estimate in \refeq{ind}, and
then invoke \refprop{Poissonsum} to handle each of the $m-1$ new sums
that arise.  Thus, the contribution of the main term above is, for
some $|\theta|\le 1$,
 \begin{displaymath}
   \begin{array}{r}
     \displaystyle{%
       \frac{r_1!\cdots r_m!}{(r_1+\ldots+r_m+m-1)!} 
       (a+\xi_1+\ldots+\xi_m)^{r_1+\ldots+r_m+m-1} +} \\
       \displaystyle{%
       \frac{\theta}{3(r-1)}  r_1!\cdots r_m!
       \frac{(a+\xi_1+\ldots+\xi_m)^{r_1+\ldots+r_m+m-3}}{(r_1+\ldots+r_m+m-3)!}},
   \end{array}
 \end{displaymath}   
 while the $j$-th term on the RHS contributes 
 \begin{displaymath}
   \begin{array}{r}
     \displaystyle{r_1! \cdots r_m! \binom{m-2}{j}
       \Big(\frac{1}{3(r-1)}\Big)^j
       \Big(\frac{(a+\xi_1+\ldots+\xi_m)^{r_1+\ldots+r_m+m-1-2j}}{(r_1
         + \ldots+r_m+m-1-2j)!}} \\
     \displaystyle{%
       + \frac{1}{3(r-1)} \frac{(a+\xi_1+\ldots+
         \xi_m)^{r_1+\ldots+r_m-1-2j-2}}{(r_1+\ldots+r_m+m-1-2j-2)!} \Big)}.
   \end{array}
 \end{displaymath}
 Using these in \refeq{ind} and the above estimate, and using the
 triangle inequality, and that $\binom{m-1}{j}= \binom{m-2}{j}+
 \binom{m-2}{j-1}$ we obtain the Proposition.
\end{proof}
 
Consider now a general deck of $n$ cards with $D_1$ $1$'s followed by
$D_2$ $2$'s and so on ending with $D_m$ $m$'s.  Recall that the
separation is maximum for the reverse configuration of the deck, and
that probability is given in \refthm{sepG}.  We now use \refprop{GPS}
to find asymptotics for that separation distance.  The following is
our `rule of thumb.'
 
\begin{theorem} 
  Consider a deck of $n$ cards of $m$-types as above.  Suppose that
  $D_i\ge d \ge 3$ for all $1\le i \le m$.  Then the separation
  distance is
  \[
  1- (1+\eta) \frac{a^{m-1}}{(n+1)\cdots(n+m-1)} \sum_{j=0}^{m-1} (-1)^j
  \binom{m-1}{j} \Big(1-\frac ja\Big)^{n+m-1} ,
  \]
  where $\eta$ is a real number satisfying 
  \[ 
  |\eta| \le \Big(1 +\frac{n^2}{3(d-2)(a-m+1)^2}\Big)^{m-1} -1.
  \]
  \label{thm:ROT}  
\end{theorem} 

\begin{proof}   
  Recall the expression for the separation distance given in
  \refthm{sepG}.  To evaluate this, we require an understanding of 
  \begin{displaymath}
    \begin{array}{l}
      \displaystyle{%
      \sum_{\substack{a_1+ \ldots + a_m =a \\ a_j \ge 1}}
      a_m^{D_m} \prod_{j=1}^{m-1} (a_j^{D_j}-(a_j-1)^{D_j})}   \\
      \hspace{3em} \displaystyle{%
      = \int_0^1 \cdots \int_0^1 \sum_{\substack{a_1
          +\ldots+a_m=a \\ a_j\ge 1}} a_m^{D_m} \prod_{j=1}^{m-1}
      \Big( D_j (a_j-1+\xi_j)^{D_j-1} d\xi_j \Big)}.
    \end{array}
  \end{displaymath}    
  We now invoke \refprop{GPS}.  Thus the above equals for some
  $|\theta| \le 1$
  \begin{displaymath}
    \begin{array}{l}
      \displaystyle{%
        \prod_{j=1}^{m} D_j!  \int_0^1 \cdots \int_0^1 \Big( 
        \frac{(a-(m-1)+\xi_1+\ldots+\xi_{m-1})^{n}}{n!} + }\\
      \displaystyle{+ \theta \sum_{j=1}^{m-1} \!\!
        \binom{m\!-\!1}{j} \! \Big(\frac{1}{3(d\!-\!2)}\Big)^j
        \frac{(a\!-\!(m\!-\!1)\!+\!\xi_1\!+\!\ldots\!+\!
          \xi_{m-1})^{n-2j}}{(n-2j)!}\Big) d\xi_1 \!\cdots\! d\xi_{m-1}.}
    \end{array}
  \end{displaymath}
  We may simplify the above as
  \begin{eqnarray*} 
    &&\Big( 1+ \theta \Big\{ \Big(1+\frac{n^2}{3(d-2)(a-m+1)^2}\Big)^{m-1}-1\Big\} 
    \Big) \frac{D_1! \cdots D_m!}{n!}\\
    &\times& \int_0^1\ldots\int_0^1 (a-m+1 +\xi_1+\ldots+\xi_{m-1})^{n} 
    d\xi_1 \cdots d\xi_{m-1},\\
  \end{eqnarray*}
and evaluating the integrals above this is 
\begin{displaymath}
  \begin{array}{r}
    \displaystyle{%
      \left( 1+ \theta \Big\{\left( 1+\frac{n^2}{3(d-2)(a-m+1)^2}
        \right)^{m-1}-1\Big\} \right) \cdot} \\
    \displaystyle{%
      \frac{D_1! \cdots D_m!}{n!} \sum_{j=0}^{m-1} (-1)^j
      \binom{m-1}{j} (a-j)^{n-m+1}}.
  \end{array}
\end{displaymath}
The Theorem follows.
\end{proof}
 
 \begin{remark}  For simplicity we have restricted ourselves to 
the case when each pile has at least three cards.  With more 
effort we could extend the analysis to include doubleton piles.  
The case of some singleton piles needs some modifications 
to our formula, but this variant can also be worked out. 
\end{remark}

 \begin{remark}  From \refthm{ROT} one can show that 
for a general decks as above, one needs $a$ of size about $nm$ before 
the separation distance becomes small.  We note that when $a$ is 
of size about $nm$, the quantity $\eta$ appearing in \refthm{ROT} 
is of size about $1/(m(d-2))$, so that the estimates furnished 
above represent a true asymptotic unless both $m$ and $d$ happen 
to be small.  In other words, when we either have many piles, 
or a small number of thick piles, \refthm{ROT} gives a good asymptotic.
\end{remark}

\begin{remark}
\label{rmk:ROT}
  While asymptotic, \refthm{ROT} is astonishingly accurate for decks
  of practical interest. For example, comparing exact calculations in
  \reftab{sep} with approximations using this rule of thumb in
  \reftab{thumb} shows that after only $3$ shuffles, the numbers agree
  to the given precision. Moreover, the simplicity of the formula in
  \refthm{ROT} allows much further computations than are possible
  using the formula in \refthm{sepG}.

  We now give a heuristic for why our rule of thumb 
is numerically so accurate; this was hinted at previously in our 
remark following \refcor{sepnn-bound}.  Let $k\ge  0$ be an integer, and define 
\[
f_k(z) = \sum_{r=0}^{\infty} r^{k} z^r,
\] 
with the convention that $0^0=1$.  Thus $f_0(z) =1/(1-z)$, $f_1(z) =
z/(1-z)^2$, and in general $f_k(z) = A_k(z)/(1-z)^{k+1}$ where
$A_k(z)$ denotes the $k$-th Eulerian polynomial.  The sum over $a_1$,
$\ldots$, $a_m$ appearing in our proof of \refthm{ROT} is simply the
coefficient of $z^{a}$ in the generating function
$(1-z)^{m-1}f_{D_1}(z)\cdots f_{D_m}(z)$.  Our rule of thumb may be
interpreted as saying that
\begin{equation}
(1-z)^{m-1} f_{D_1}(z)\cdots f_{D_m}(z) \approx 
\frac{D_1! \cdots D_{m}!}{(n+m-1)!} (1-z)^{m-1} f_{n+m-1}(z).
\label{eqn:genROT}
\end{equation}
To explain the sense in which \refeq{genROT} holds, note
that $f_k(z)$ extends meromorphically to the 
complex plane, and it has a pole of order $k+1$ at $z=1$.  
Moreover it is easy to see that $f_k(z) - k!/(1-z)^{k+1}$ 
has a pole of order at most $k$ at $z=1$.  Therefore, 
the LHS and RHS of \refeq{genROT} have 
poles of order $n+1$ at $z=1$, and their leading order contributions 
match.  Therefore the difference between the RHS and LHS 
of \refeq{genROT} has a pole of order at most $n$ at $z=1$.  
But in fact, this difference can have a pole of order at 
most $n-d$ at $z=1$, and thus the approximation in \refeq{genROT} 
is tighter than what may be expected {\sl a priori}.  To obtain 
our result on the order of the pole, we record that 
one can show 
\[ 
f_k(z) = \frac{k!}{(1-z)^{k+1}} \Big( \frac{(z-1)}{ \log z} \Big)^{k+1} +\zeta(-k) + 
O(1-z).
\]
\end{remark}

\section{Comparing $2$-shuffles with different starting patterns}
\label{sec:alternating}

Conger and Viswanath note that the initial configuration can affect
the speed of convergence to stationary. In this section, we
investigate this for a deck with $n$ red and $n$ black cards. Consider
first starting with reds on top. If the initial cut is at $n$ (the
most likely value) then the red-black pattern is perfectly mixed after
a single shuffle. More generally, by \refcor{2shuffle}, the chance of
the deck $w$ resulting from a single $2$-shuffle of a deck with $n$
red cards atop $n$ black cards is given by
\begin{displaymath}
    Q_2(w) = \frac{1}{2^{2n}} \left( 2^{\head(w)} + 2^{\tail(w)} - 1 \right).
\end{displaymath}

Consider next the result of $2$-shuffles on the {\em alternating deck}
red-black-red-black-$\cdots$. As motivation, we recall a popular card
trick: Begin with a deck of $2n$ cards arranged alternately red,
black, red, black, etc. The deck may be cut any number of times. Have
the deck turned face up and cut (with cuts completed) until one of the
cuts results in the two piles having cards of opposite color
uppermost. At this point, ask one of the participants to riffle
shuffle the two piles together. The resulting arrangement has the top
two cards containing one red and one black, the next two cards
containing one red and one black, and so on throughout the deck. This
trick is called the Gilbreath Principle after its inventor, the
mathematician Norman Gilbreath. It is developed, with many variations,
in Chapter 4 of \cite{Gardner1966}. 

From the trick we see that beginning with an alternating deck severely
limits the possibilities. Which start mixes faster? The following
developments both explain the trick and give a useful formula for
analysis.

\begin{lemma}
  The number of deck patterns resulting from a cut with an odd number
  of cards in both piles followed by a riffle shuffle is
  $2^{n}$. Similarly, the number of deck patterns resulting from a cut
  with both piles even followed by a riffle shuffle is $2^{n-1}$.
\label{lem:drops}
\end{lemma}

\begin{proof}
  For the case of an odd cut, the last two cards after the riffle
  shuffle must be a red and a black card. No matter what piles these
  two cards fell from, the next two cards must also consist of one red
  and one black card. Continuing on, the possible resulting decks are
  exactly those where the $i$th and $i+1$st cards have different
  colors for $i=1,3,\ldots,2n-1$. The number of such decks is exactly
  $2^n$, since each of the order of each of the $n$ pairs is
  independent. 

  For an even cut, we proceed by induction noting that the case when
  $n=1,2,3$ are easily solved by inspection. In this case, the only
  resulting decks will necessarily begin with a red card and end with
  a black card. The number of decks beginning with two red cards or
  ending with two black cards is determined by the previous case since
  removing the top or bottom card from each pile results in piles with
  an odd number of cards, giving $2^{n-1}$ possibilities. However, we
  must discount the over counted case of decks beginning with two red
  cards and ending with two black cards, and, by induction since the
  piles are again both even, there are $2^{n-3}$ such decks. Finally,
  the remaining case must be decks beginning and ending with a red
  card followed by a black card. In this case, again, the piles remain
  even and by induction the number of such decks is
  $2^{n-3}$. Therefore the total count for cuts with both piles even
  is $2^{n-1} - 2^{n-3} + 2^{n-3} = 2^{n-1}$.
\end{proof}

The proof of the lemma shows exactly why the card trick is a success:
to have different colors on the top of the two piles, the cut must
have been odd. Therefore the first two cards dropped consist of one
red and one black, and the next two cards dropped consist of one red
and one black, and so on. Also from the lemma, we see that the only
deck that can result from either an odd cut or an even cut is the
identity.

\begin{proposition}
  The chance of a $2$-shuffle of the alternating deck resulting in a
  deck configuration $w$ is given by
  \begin{equation}
    2^{2n} \cdot Q_2(w) = \left\{ \begin{array}{cl}
        2^{n-1} + 2^n & \mbox{if} \ w = w_0 \\
        2^{n-1}       & \mbox{if} \ w \in O\setminus w_0, \\
        2^n          & \mbox{if} \ w \in E\setminus w_0, \\
        0            & \mbox{otherwise},
      \end{array} \right.
  \label{eqn:Q2-2}
  \end{equation}
  where $w_0$ is the initial alternating deck and $O$ (respectively,
  $E$) is the set of decks that can result from riffling together the
  two piles from cutting the alternating deck when both piles have an
  odd (respectively, even) number of cards.
\label{prop:RBRB}
\end{proposition}

\begin{proof}
  Let $w,u \in O$. Then the total number of ways $w$ can result from
  any odd cut is equal to the total number of ways $u$ can result from
  any odd cut. The same is true replacing $O$ with $E$ and ``odd''
  with ``even''. From the binomial identity 
  \begin{displaymath}
    \sum_{k=0}^{n} (-1)^k \binom{n}{k} = 0 \hspace{2ex} \leadsto \hspace{2ex}
    \sum_{k \ \mathrm{odd}} \binom{2n}{k} = \sum_{k \ \mathrm{even}}
    \binom{2n}{k},
  \end{displaymath}
  we must have both the right-hand sums equal to $2^{2n-1}$. Therefore,
  by \reflem{drops}, the total number of ways $w$ can result from an
  odd cut (assuming it can) is $2^{2n-1}/2^{n} = 2^{n-1}$, and,
  similarly, the total number of ways $w$ can result from an even cut
  (assuming it can) is $2^{2n-1}/2^{n-1} = 2^{n}$.
\end{proof}

It follows from \refeq{Q2-2} that the separation distance for a
$2$-shuffle is $\SEP(2)=1$ when $n \geq 3$. Furthermore, since
$\binom{2n}{n} \geq 2^{n}$, we can compute the total variation of a
$2$-shuffle to be
\begin{equation}
  \left\|Q_2 - U\right\|_{TV} = \frac{1}{2} \left( 1 - \frac{2^n +
      2^{n-1} - 1}{\binom{2n}{n}} \right),
  \label{eqn:TValt}
\end{equation}
which goes to $.5$ exponentially fast as $n$ goes to infinity. In
contrast, starting with reds above blacks, asymptotic analysis of
\refeq{2sh} shows that the total variation tends to $1$ after a single
shuffle when $n$ is large. Thus an alternating start leads to faster
mixing.

\bigskip
\begin{center}
{\sc Acknowledgements}
\end{center}

The authors thank Jason Fulman for careful comments and references.
We thank Mark Conger and Divakar Viswanath for extensive, useful
correspondence and for sending us copies of \cite{CoVi,Conger2007}
while our paper was under review, and Steve Lalley for help beyond the
call of duty. We also thank an extremely helpful referee for detailed
comments and suggestions. Finally, we thank MSRI and the participants
of the combinatorial representation theory program where this work
began.

  \bibliographystyle{abbrv} 
  \bibliography{../references}


%
%
\appendix

\section{Random walks on groups}
\label{sec:groups}

In this appendix, we reformulate shuffling in terms of random walks
on the symmetric group $\Sn$, so that our investigation of particular
properties of a deck becomes the quotient walk on Young subgroups of
$\Sn$. 

Let $G$ be a finite group with $Q(g) \geq 0$, $\sum_{g \in G} Q(g) =
1$ a probability on $G$. The walk in \refeq{convolution} may be called
the {\em left walk} since it consists of repeatedly picking elements
independently with probability $Q$, say $g_1, g_2, g_3, \ldots$, and,
starting at the identity $1_{G}$, multiplying on the left by
$g_i$. The generates a {\em random walk on $G$},
\begin{displaymath}
  1_G, \ g_1, \ g_2 g_1, \ g_3 g_2 g_1, \ \ldots .
\end{displaymath}
By inspection, the chance that the walk is at $g$ after $k$ steps is
$Q^{*k}(g)$, where $Q^0(g) = \delta_{1_G,g}$.

An algebraic method of focusing on aspects of the walk is to use the
{\em quotient walk}. Let $H \leq G$ be a subgroup of $G$, and set $X =
G/H = \{xH\}$ to be the set of left cosets of $H$ in $G$. The quotient
walk is derived from the walk above by simply reporting to which coset
the current position of the walk belongs. The quotient walk is a
Markov chain on $X$ with transition matrix given by
\begin{equation}
  K(x,y) = Q(y H x^{-1}) = \sum_{h \in H} Q(yhx^{-1}) .
\label{eqn:transition}
\end{equation}
Note that $K$ is well-defined (i.e. independent of the choice of coset
representatives) and that $K$ is doubly stochastic. Thus the uniform
distribution on $X$, $U(x) = |H|/|G|$, is a stationary distribution
for $K$. The chain $K$ is reversible if and only if $Q$ is symmetric
(i.e. $Q(g) = Q(g^{-1})$). Note that this is not the case for riffle
shuffles. While intuitively obvious, the following shows the basic
fact that powers of the matrix $K$ correspond to convolving and taking
cosets.

\begin{proposition}
  For $Q$ a probability distribution on a finite group $G$ and $K$ as
  defined in \refeq{transition}, we have
  \begin{displaymath}
    K^{l}(x,y) = Q^{*l} (yHx^{-1}) .
  \end{displaymath}
\label{prop:convolve}
\end{proposition}

\begin{proof}
  The result is immediate from the definitions for $l=0,1$. We prove
  the result for $l=2$, the general case being similar. Note that
  \begin{displaymath}
    K^2(x,y) = \sum_z K(x,z) K(z,y) = \sum_z \sum_{h_1, h_2} Q(z h_1
    x^{-1}) Q(y h_2 z^{-1}) .
  \end{displaymath}
  Setting $h_2 = h h_1^{-1}$, noting that $zh_1$ runs over $G$ as $z$
  runs over $X$ and $h_1$ over $H$, and setting $g_1 = gx^{-1}$, we
  have
  \begin{eqnarray*}
    K^2 (x,y) & = & \sum_h \sum_g Q(gx^{-1}) Q(y h g^{-1}) \\
    & = & \sum_h \sum_{g_1} Q(g_1) Q(yh x^{-1} g_{1}^{-1}) = Q^2(yHx^{-1}).
  \end{eqnarray*}
\end{proof}

We may identify permutations in $\Sn$ with arrangements of a deck of
$n$ cards by setting $\sigma(i)$ to be the label of the card at
position $i$ from the top. Thus the permutation $2 \ 1 \ 4 \ 3$ is
associated with four cards where ``2'' is on top, followed by ``1'',
followed by ``4'', and finally ``3'' is on the bottom. If we consider
the cards labelled $1,2,\ldots,k$ to be ``red'' cards, and the cards
labelled $k+1,k+2,\ldots,n$ to be ``black'' cards, with all cards of the
same color indistinguishable, the coset space
\begin{displaymath}
  X = \Sn / (\mathcal{S}_{k} \times \mathcal{S}_{n-k})
\end{displaymath}
is naturally associated with the $\binom{n}{k}$ arrangements of red
and black {\em unlabeled} cards. Here, of course, we identify an element
of $\mathcal{S}_{k} \times \mathcal{S}_{n-k} \leq \Sn$ as permuting
the first $k$ and last $n-k$ cards among themselves. Similar
constructions work for suits or values. Thus \refprop{convolve} shows
that the processes studied in the body of this paper are Markov
chains.

\section{Shuffling by random transpositions}
\label{app:reps}

Let $L^2(X) = \{f : X \rightarrow \mathbb{C} \}$ be the set of
complex-valued functions on $X$ with inner product defined by
\begin{equation}
  \langle f_1 | f_2 \rangle = \frac{1}{|X|} \sum_x f_1(x)
  \overline{f_2(x)} . 
\label{eqn:inner}
\end{equation}
If $K$ is symmetric, then real-valued  functions may be used. The
transition matrix $K$ operates on $L^2$ via 
\begin{equation}
  K f(x) = \sum_y K(x,y) f(y) .
\label{eqn:Kf}
\end{equation}
In the present case, $L^2(X) = \Ind_{H}^{G}(\mathbf{1})$, the usual
permutation representation of $G$ acting on left cosets $X = G/H$,
with $T_g f(x) = f(g^{-1}x)$. By construction, the action of $G$
commutes with $K$, i.e.
\begin{equation}
  T_{g} (Kf) = K(T_{g} f)
\label{eqn:KGcommute}
\end{equation}
for all $f \in L^2(X)$ and all $g \in G$. This implies that group
representation theory can be used to reduce the operator $K$ (or
diagonalize $K$ in the case when $K$ is symmetric). This classical
topic is well developed in F{\"a}ssler-Steifel \cite{FaSt1992} and
Boyd, et. al. \cite{BD2005}.

Let $\hG$ denote the set of irreducible representations of the finite
group $G$. For $\rho \in \hG$, the Fourier transform of $f \in L^2(G)$
at $\rho$ is defined by
\begin{displaymath}
  \widehat{f}(\rho) = \sum_{g \in G} f(g) \rho(g) .
\end{displaymath}
As usual, Fourier transform turns convolution into products, i.e.
\begin{displaymath}
  \widehat{Q^{*k}}(\rho) = \widehat{Q}(\rho)^k.
\end{displaymath}
Schur's lemma implies that the uniform distribution has zero transform
\begin{displaymath}
  \widehat{U}(\rho) = \left\{ \begin{array}{ll}
      1 & \mbox{if $\rho$ is trivial,} \\
      0 & \mbox{otherwise.} \\
    \end{array} \right.
\end{displaymath}
The Fourier inversion theorem reconstructs $f$ from
$\{\widehat{f}(\rho)\}$ by
\begin{displaymath}
  f(g) = \frac{1}{|G|} \sum_{\rho \in \hG} \mathrm{dim}_{\rho}
  \mathrm{Tr}\left( \widehat{f}(\rho) \rho(g^{-1}) \right) .
\end{displaymath}
For background, see Serre \cite{Serre1977}, Diaconis
\cite{Diaconis1988} or Ceccherini, et. al \cite{CST2008} where
many applications are given.

Suppose the induced representation $L^2(X)$ decomposes into
irreducibles as 
\begin{equation}
  L^2(X) = \bigoplus_{\rho \in \hG} V_{\rho}^{\oplus a_{\rho}} .
\label{eqn:L2V}
\end{equation}
Then since $K$ commutes with $G$, $K$ sends each of the spaces
$V_{\rho}^{\oplus a_{\rho}}$ into itself. Further reductions may be
possible if $Q$ has suitable symmetries. The following widely studied
special case is relevant.

\begin{definition}
  The pair $H \leq G$ is a {\em Gelfand pair} if $L^2(X)$ is
  multiplicity free, i.e. all $a_{\rho}$ in \refeq{L2V} are either $0$
  or $1$.
\label{defn:gelfand}
\end{definition}

For example, when $1 \leq k \leq n/2$, $\mathcal{S}_{k} \times
\mathcal{S}_{n-k} \leq \Sn$ is a Gelfand pair with
\begin{equation}
  L^2(X) = \bigoplus_{i=0}^{k} S^{n-i,i}.
\label{eqn:Skn-k}
\end{equation}
Recall that the irreducible representations of $\Sn$ are indexed by
partitions $\lambda$ of $n$. If $S^{\lambda}$ denotes the $\lambda$th
representation (Specht modules), the sum in \refeq{Skn-k} runs over
partitions into two parts with the smaller part at most $k$. For
further background on Gelfand pairs, including examples and
applications, see \cite{Diaconis1988,CST2008}.

Now we study a deck of red and black cards after repeated random
transposition shuffles. Recall that Diaconis-Shahshahani
\cite{DiSh1981} show that it takes $\frac{1}{2}n(\log(n)+c)$
shuffles to mix $n$ distinct cards. To be precise, the measure on
$\Sn$ that drives the walks is
\begin{displaymath}
  Q(\sigma) = \left\{ \begin{array}{cl}
      1/n   & \mbox{if} \ \sigma = \mathrm{id}, \\
      2/n^2 & \mbox{if} \ \sigma = (i,j), \\
      0     & \mbox{otherwise.}
    \end{array} \right.  
\end{displaymath}
Throughout the following, all walks begin at the identity permutation,
and we use the convention that $\pi(i)$ is the label of the card at
position $i$.

First, we follow the position of the top card; i.e. the two of hearts
is the only red card followed by $n-1$ black cards. The transition
matrix for this walk is given by
\begin{equation}
  P(i,j) = \left\{ \begin{array}{cl}
    \displaystyle{\frac{1}{n} + \frac{(n-2)(n-3)}{n^2}} 
    & \mbox{if} \ i=j, \\[\vsp]
    \displaystyle{\frac{2}{n^2}} & \mbox{if} \ i \neq j.
    \end{array} \right.
\label{eqn:Pij}
\end{equation}
Note that this is symmetric, with $\Pi(i) = 1/n$ as the stationary
distribution. 

\begin{proposition}
  For the transition matrix $P(i,j)$ above and all $l \geq 0$, we have 
  \begin{equation}
    P^l(i,j) = \left\{ \begin{array}{cl}
        \displaystyle{\frac{1}{n} + \left(1 - \frac{2}{n}\right)^l
          \left(1 - \frac{1}{n}\right)} & \mbox{if} \ i=j, \\[\vsp]
        \displaystyle{\frac{1}{n} - \left(1 - \frac{2}{n}\right)^l
          \frac{1}{n}} & \mbox{if} \ i \neq j.
      \end{array} \right.
    \label{eqn:Plij}
  \end{equation}
  From this it follows that 
  \begin{displaymath}
    \SEP(l) = \left(1 - \frac{2}{n} \right)^l \hspace{2ex} \mbox{and}
    \hspace{2ex} \| P - \Pi \|_{TV} = \left(1 - \frac{2}{n} \right)^l
    \left(1 - \frac{1}{n}\right) .
  \end{displaymath}
  \label{prop:Pij}
\end{proposition}

\begin{proof}
  The results for the separation and total variation distances follow
  from \refeq{Plij} and the definitions. It is possible to give a
  direct combinatorial argument for \refeq{Plij}, but the following
  representation theoretic argument generalizes readily to find
  similar formula for $j$-tuples of cards.

  The random transposition measure $Q$ is constant on conjugacy
  classes of $\Sn$ and so acts on each irreducible representation as a
  constant times the identity. These constants are given explicitly by
  Diaconis-Shahshahani \cite{DiSh1981}, involving characters and
  dimensions of the representation. Consider the operator
  $K(\sigma,\tau) = Q(\tau\sigma^{-1})$ on the regular
  representation. The function $f(\sigma) = \delta_{1,\sigma(i)} -
  1/n$ lies in the $n-1$ copies of the $n-1$-dimensional
  representation corresponding to the partition $(n-1,1)$. The
  operator $K$ acts on this space by multiplication by $1 - 2/n$. Thus
  \begin{eqnarray*}
    P_{\sigma} \left( \begin{array}{c}
        \mbox{card labelled 1} \\ \mbox{at position $i$} \\ \mbox{after
          $l$ shuffles} \end{array} \right) - \frac{1}{n} & = &
    K^{l} f(\sigma) = \left(1 - \frac{2}{n} \right)^{l} f(\sigma) \\
    & = & \left(1 - \frac{2}{n} \right)^{l} \left( \delta_{1,\sigma(i)} -
      \frac{1}{n} \right).
  \end{eqnarray*}
  Here $\sigma$ is the starting arrangement. Evaluating the right-hand
  side gives \refeq{Plij}.
\end{proof}

Next we consider the deck with $N=2n$ cards where the (original) top
$n$ cards are red and the (original) bottom $n$ cards are black. In
this case, we think of the the random transposition operator acting on
the quotient space $\mathcal{S}_{N} / \Sn \times \Sn$. For $x,y \in
\mathcal{S}_{N} / \Sn \times \Sn$, the induced Markov chain is
\begin{equation}
  K(x,y) = \left\{ \begin{array}{cl}
    \frac{1}{N^2} & \mbox{if $x \neq y$ differ by a transposition,} \\
    \frac{1}{N} + \frac{(n(n-1))^2}{N^2} & \mbox{if $x=y$,} \\
    0 & \mbox{otherwise}.
  \end{array} \right.
\label{eqn:Kxy}  
\end{equation}
This chain has uniform stationary distribution $\Pi(x) =
1/\binom{N}{n}$.

The chain $K$ is invariant under $\mathcal{S}_N$, i.e. $K(x,y) =
K(\sigma x,\sigma y)$, so the distance to stationary does not depend
on the original configuration. As noted earlier, the pair $\Sn \times
\Sn, \mathcal{S}_N$ is a Gelfand pair, so \refeq{Skn-k} allows an easy
determination of the eigen values and rate of convergence. 

\begin{proposition}
  For the Markov chain $K$ on $\mathcal{S}_{N} / \Sn \times \Sn$, the
  eigen values are
  \begin{displaymath}
    \beta_0 = 1, \hspace{2ex} \beta_j = \frac{1}{N} + \frac{1}{N^2}
    \left((N-j)^2 - (N-j) + j^2 -3j \right),
  \end{displaymath}
  $j = 1, \ldots, n$. The multiplicity of $\beta_j$ is $m_j =
  \binom{N-1}{j}$. Moreover, there is a universal constant $A$ such
  that if $l = \frac{1}{4}N(\log N + C)$, then
  \begin{displaymath}
    \left\| K^l - \Pi \right\|_{TV} \leq A e^{-c/2}.
  \end{displaymath}
\label{prop:GelRB}
\end{proposition}

\begin{proof}
  The operator $K$ acts on $L^2 \left(\mathcal{S}_{N} / \Sn \times \Sn
  \right)$ as the element of the group algebra
  \begin{displaymath}
    \frac{1}{N} \mathrm{Id} + \frac{2}{N^2} \sum_{i<j}(i,j).
  \end{displaymath}
  As shown in \cite{DiHo2002}, this element acts on the irreducibles
  $\mathcal{S}^{n-j,j}$ as a constant times the identity, with the
  constant being $\beta_j$ and the multiplicity being the dimension of
  $\mathcal{S}^{n-j,j}$. This proves the first part.

  The remaining claims can be proved following the argument in
  \cite{DiHo2002}: bound the total variation distance by the $L^2$
  norm, express this in terms of the eigen values and average over the
  starting state. This reduces the problem to bounding
  \begin{displaymath}
    \sum_{j=1}^{n} m_j \beta_j^{2l}.
  \end{displaymath}
  The lead term in this is 
  \begin{displaymath}
    (N-1) \left( 1 - \frac{2}{N} \right)^{2l} \leq e^{-c}.
  \end{displaymath}
  For $l$ of the form $\frac{1}{4}N(\log N + c)$, the other terms
  are smaller and sum in a reasonably standard fashion. The terms are
  the same as in \cite{DiHo2002}, so we suppress further details.
\end{proof}

\begin{remark}
  It is easy to give a lower bound showing that after
  $l=\frac{1}{4}N(\log N + c)$ steps the distance to stationary is
  bounded away from $0$ for large $N$. Further, in this case, the
  distance tends to $1$ if $c = c_N$ tends to $-\infty$.
\end{remark}

These results show that for red-black mixing, there is a total
variation cutoff at $\frac{1}{4}N\log N$. Note that single card mixing
does not have a cutoff, recalling that in \refprop{Pij} the deck has
size $n$ and in \refprop{GelRB} the deck has size $N=2n$.

\end{document}